\documentclass{amsart}

\usepackage{epsfig, amssymb, latexsym, graphicx, color}
\input xy

\xyoption{all}

\theoremstyle{plain}
\newtheorem{theorem}{Theorem}
\newtheorem*{maintheorem}{Density Theorem}
\newtheorem{lemma}[theorem]{Lemma}
\newtheorem{proposition}[theorem]{Proposition}
\newtheorem{corollary}[theorem]{Corollary}
\newtheorem{definition}[theorem]{Definition}

\theoremstyle{definition}
\newtheorem{examples}{Examples}
\newtheorem{example}[examples]{Example}

\newcommand{\s}{\ensuremath{\sigma}}

\newcommand{\Z}{\ensuremath{\mathbb{Z}}}

\newcommand{\Hyp}{\ensuremath{\mathbb{H}}}

\newcommand{\cC}{\ensuremath{\mathcal{C}}}
\newcommand{\cS}{\ensuremath{\mathcal{S}}}

\newcommand{\cK}{\ensuremath{\mathcal{K}}}

\newcommand{\Aut}{\ensuremath{\operatorname{Aut}}}

\newcommand{\Ad}{\ensuremath{\operatorname{Ad}}}

\newcommand{\Lk}{\ensuremath{\operatorname{Lk}}}

\newcommand{\Ch}{\ensuremath{\operatorname{Ch}}}

\newcommand{\Stab}{\ensuremath{\operatorname{Stab}}}
\newcommand{\CommG}{\mathrm{Comm}_G}
\newcommand{\Comm}{\ensuremath{\operatorname{Comm}}}
\newcommand{\CAT}{\ensuremath{\operatorname{CAT}}}

\newcommand{\Id}{\ensuremath{\operatorname{Id}}}
\newcommand{\Sym}{\ensuremath{\operatorname{Sym}}}

\newcommand{\G}{\Gamma}
\newcommand{\bC}{\partial \cC} 

\newcommand{\bs}{\backslash}

\def\polhk#1{\setbox0=\hbox{#1}{\ooalign{\hidewidth
    \lower1.0ex\hbox{$\,\lhook$}\hidewidth\crcr\unhbox0}}}
\newcommand{\Swiatkowski}{\'Swi{\polhk{a}}tkowski}

\begin{document}

\title[Density of commensurators]{Density of commensurators for uniform lattices of right-angled buildings}


\author{Angela Kubena}
\address{Department of Mathematics, University of Michigan,
2074 East Hall, 530 Church Street, Ann Arbor, MI 48109-1043, USA}
\email{akubena@umich.edu}\thanks{The second author was supported in part by NSF Grant No. DMS-0805206 and is now supported in part by ARC Grant No. DP110100440.}

\author{Anne Thomas}
\address{School of Mathematics and Statistics, Carslaw Building F07, University of Sydney NSW 2006, Australia}
\email{anne.thomas@sydney.edu.au}

\begin{abstract} Let $G$ be the automorphism group of a regular right-angled building $X$.  The
``standard uniform lattice" $\G_0 \leq G$ is a canonical graph product of finite groups, which acts
discretely on $X$ with quotient a chamber.  We prove that the commensurator of $\Gamma_0$ is
dense in $G$.  This result was also obtained by Haglund \cite{H_FI}.  For our proof, we develop carefully a technique of ``unfoldings" of complexes of groups.  We use unfoldings to construct a sequence of uniform lattices $\G_n \leq G$, each commensurable to
$\G_0$, and then apply the theory of group actions on complexes of groups to the sequence
$\G_n$.  As further applications of unfoldings, we determine exactly when the group $G$ is
nondiscrete, and prove that $G$ acts strongly transitively on $X$. \end{abstract}

\maketitle

\section*{Introduction}\label{s:intro}

Two subgroups $\G_0$ and $\G_1$ of a group $G$ are \emph{commensurable} if the intersection $\G_0 \cap
\G_1$ has finite index in both $\G_0$ and $\G_1$.  The \emph{commensurator} of $\G \leq G$ in $G$ is the
group \[\CommG(\G) := \{ g \in G \mid \mbox{$g\G g^{-1}$ and $\G$ are commensurable} \}. \]  Note that
$\Comm_G(\G)$ contains the normalizer $N_G(\G)$.  It is a classical fact that if $G$ is a connected
semisimple Lie group, with trivial center and no compact factors, and $\G \leq G$ is an irreducible
lattice, then either $\G$ is finite index in $\CommG(\G)$ or $\CommG(\G)$ is dense in $G$
(see~\cite{Z}).  Moreover Margulis~\cite{M} proved that $\G$ is arithmetic if and only if $\CommG(\G)$
is dense.

A semisimple Lie group is a locally compact topological group.  If $X$ is a locally finite, simply
connected polyhedral complex, then the group $G = \Aut(X)$ is also locally compact.  It turns out that a subgroup
$\G \leq G$ is a uniform lattice in $G$  if and only if $\G$ acts cocompactly on $X$ with finite cell
stabilizers (see Section~\ref{ss:lattices}). Lattices in such groups $G$ share many properties with
lattices in semisimple Lie groups, but also exhibit new and unexpected phenomena (see the
surveys~\cite{Lu} and~\cite{FHT}).

In this setting, the one-dimensional case is $X$ a locally finite tree.  Liu~\cite{Liu} proved that
the commensurator of the ``standard uniform lattice" $\G_0$ is dense in $G=\Aut(X)$; here $\G_0$ is a canonical
graph of finite cyclic groups over the finite  quotient $G \bs X$.  In addition,
Leighton~\cite{Leighton} and Bass--Kulkarni~\cite{BK} proved that all uniform lattices in $G$ are
commensurable (up to conjugacy).  Hence all uniform tree lattices have dense commensurators.  In
dimension two, Haglund~\cite{H_RCD} showed that for certain $2$--dimensional Davis complexes $X=X_W$,
the Coxeter group $W$, which may be regarded as a uniform lattice in $G=\Aut(X)$, has dense
commensurator.

We consider higher-dimensional cases, focusing on regular right-angled buildings $X$ (see
Section~\ref{ss:rabs}).  Such buildings exist in arbitrary dimension.  Examples include products of
finitely many regular trees, and Bourdon's building $I_{p,q}$, the unique $2$--complex in which every
$2$--cell is a regular right-angled hyperbolic $p$--gon and the link of each vertex is the complete
bipartite graph $K_{q,q}$ (see~\cite{B}).  The ``standard uniform lattice" $\G_0 \leq G = \Aut(X)$,
defined in Section~\ref{ss:complexes_of_groups} below, is a canonical graph product of finite cyclic
groups, which acts on $X$ with fundamental domain a chamber.  Our main result is:

\begin{maintheorem} Let $G$ be the automorphism group of a locally finite regular right-angled
building $X$, and let $\G_0$ be the standard uniform lattice in $G$.  Then $\Comm_G(\G_0)$ is dense
in $G$. \end{maintheorem}

\noindent This theorem was proved independently and using different methods by Haglund~\cite[Theorem 4.30]{H_FI}.  Although the paper \cite{H_FI} was submitted in 2004, it was not publicly available and was not known to us.

In contrast to our main result, we show in Section~\ref{s:norm} below that for all $G=\Aut(X)$ with $X$ a locally finite polyhedral
complex (not necessarily a building), and all uniform lattices $\G \leq G$, the normalizer $N_G(\G)$ is
discrete.  Hence for $G$ as in the Density Theorem, the density of $\CommG(\G_0)$ does not come just
from the normalizer.

For most regular right-angled buildings $X$, it is not known whether all uniform lattices in $G=\Aut(X)$
are commensurable (up to conjugacy).  Januszkiewicz--\Swiatkowski~\cite{JS1} have established
commensurability of a class of lattices in $G$ which includes $\G_0$, where each such lattice is a graph
product of finite groups.  Hence by the Density Theorem, each such lattice has dense commensurator.  For
Bourdon's building $I_{p,q}$, Haglund~\cite{H_CS} proved that if $p \geq 6$, then all uniform lattices
in $G=\Aut(I_{p,q})$ are commensurable (up to conjugacy).  Thus by the Density Theorem, all uniform
lattices in $G$ have dense commensurators.  On the other hand, for $X$ a product of two trees,
Burger--Mozes~\cite{BM2} constructed a uniform lattice $\G \leq \Aut(X)$ which is a simple group.  It
follows that $\CommG(\G) = N_G(\G)$, which is discrete.  Thus there are cases (when $\dim(X) \geq 2$) in
which not all uniform lattices $\G \leq G = \Aut(X)$ have dense commensurators. In fact, it is an
open problem to determine whether the only possibilities for $\CommG(\G)$ are discreteness or density.
As for commensurators of nonuniform lattices in $G=\Aut(X)$, hardly anything is known, even for $X$ a tree
(see~\cite{FHT}).

If the building $X$ can be equipped with a $\CAT(-1)$ metric, then the Density
Theorem may be combined with the commensurator superrigidity theorem of Burger--Mozes~\cite{BM1} for
$\CAT(-1)$ spaces, to give rigidity results for lattices in $G=\Aut(X)$ which are commensurable to
$\G_0$.  Regular right-angled buildings with piecewise hyperbolic
$\CAT(-1)$ metrics exist in arbitrarily high dimensions~\cite{JS2}.

We now outline the proof of the Density Theorem, which is given in full in Section~\ref{s:Proof} below.
Fix a basepoint $x_0 \in X$.  Denote by $Y_n$ the combinatorial ball of radius $n$ about $x_0$ in $X$.
We first reduce the theorem to showing that for all $g \in \Stab_G(x_0)$ and for all $n \geq 0$, there
is a $\gamma_n \in \CommG(\G_0)$ such that $\gamma_n$ agrees with $g$ on the ball $Y_n$.  We then
construct a canonical uniform lattice $\G_n$ with fundamental domain the ball $Y_n$ and show that
$\G_n$ is a finite index subgroup of $\G_0$.  By considering the restriction of $g$ to $Y_n$, we are
then able to build a uniform lattice $\G_n'$ which contains a suitable element $\gamma_n$.   By our
construction, the lattice $\G_n$ is a finite index subgroup of $\G_n'$.  That is, $\G_n'$ and $\G_0$
have a common finite index subgroup $\G_n$, as sketched on the left of
Figure~\ref{f:lattice_inclusions} below.  Thus $\G_n'$ is commensurable to $\G_0$, and so $\gamma_n$
lies in $\CommG(\G_0)$, as required.

\begin{figure}[ht]
\scalebox{0.9}{\includegraphics{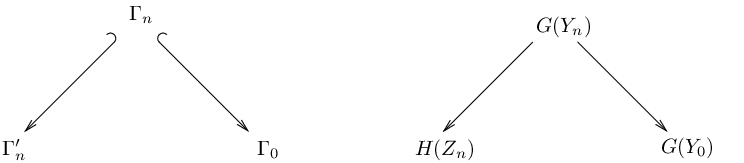}}
\caption{Inclusions of lattices (left) and coverings of complexes of groups (right)}
\label{f:lattice_inclusions}
\end{figure}

Our lattices $\G_n$ and $\G_n'$ are fundamental groups of complexes of groups (see~\cite{BH} and
Section~\ref{ss:complexes_of_groups} below).  The finite index lattice inclusions on the left of
Figure~\ref{f:lattice_inclusions} are induced by finite-sheeted coverings of complexes of groups, shown on the right of Figure~\ref{f:lattice_inclusions}.  The necessary covering theory for complexes of
groups is recalled in Section~\ref{ss:covering_theory} below (see also~\cite{BH} and~\cite{LT}).

To construct the sequence of lattices $\G_n$, in Section~\ref{s:unfolding} below we introduce and carefully develop a new tool, that of
\emph{unfoldings} of complexes of groups.  The idea of unfoldings may be viewed as a ``thicker" version of the well-known fact that convex subcomplexes of the Davis complex for a right-angled Coxeter group $W$ are fundamental domains for certain subgroups of $W$ (compare also the ``fanning" construction of Carbone~\cite{C} and the ``blocks" used in~\cite{ABJLMS}).  Now, the standard uniform lattice $\G_0$ is the fundamental
group of a complex of groups $G(Y_0)$ over a chamber $Y_0$ of $X$.  By ``unfolding" along ``sides'' of successive
unions of chambers starting from $Y_0$, and defining new local groups appropriately, we obtain a canonical family
of complexes of groups $G(Y_n)$ over the combinatorial balls $Y_n \subset X$.  The fundamental group $\G_n$ of
$G(Y_n)$ is a uniform lattice in $G = \Aut(X)$, and each $\G_n$ is a finite index subgroup of $\G_0$.  We
prove these properties of unfoldings inductively by combinatorial arguments involving close consideration of the local structure of $X$, together
with facts about
Coxeter groups, and the definition of a building as a chamber system equipped with a $W$--distance function (see
Section~\ref{ss:rabs}).

The other main tool in our proof of the Density Theorem is that of \emph{group actions on complexes of
groups}, which was introduced by the second author in~\cite{T2} (see Section~\ref{ss:group_actions}
below).  This theory is used to construct the sequence of lattices $\G_n'$, containing suitable elements
$\gamma_n$, as fundamental groups of
complexes of groups $H(Z_n)$ such that there are finite-sheeted coverings $G(Y_n) \to H(Z_n)$.

We describe in Section~\ref{ss:new_lattices} below how our two main tools, unfoldings and group actions on complexes of
groups, may be combined to construct many uniform lattices in addition to the sequences $\G_n$ and
$\G_n'$ used in the proof of the Density Theorem. To our knowledge, the lattices so obtained are new.
In particular, they do not ``come from" tree lattices, unlike the lattices in~\cite{T1}.

In Section~\ref{s:applications} below, we give two further applications of the technique of unfoldings.
First, in Theorem~\ref{t:autom_gp_intro} below, we complete the characterization of those regular right-angled buildings $X$ such that $G=\Aut(X)$ is
nondiscrete (the lattice theory of $G$ being otherwise trivial).  As we recall in
Section~\ref{ss:rabs}, the building $X$ is determined by a right-angled Coxeter system $W = \langle S
\mid (st)^{m_{st}} \rangle$ and a family of positive integers $\{q_s \}_{s \in S}$, where $q_s$ is the number of chambers of $X$ which meet at a common face of \emph{type} $s$.  (For example, in Bourdon's building $I_{p,q}$, all $q_s$ are equal to $q$.)  A polyhedral complex $L$ is said to be \emph{rigid} if for any $g \in \Aut(L)$, if $g$ fixes the star in $L$ of a vertex $\s$ of $L$, then $g$ fixes $L$.  If $L$ is not rigid, it is \emph{flexible}.  We prove:

\begin{theorem}\label{t:autom_gp_intro} Let $X$ be a regular right-angled building of type $(W,S)$ and
parameters $\{ q_s\}$. Let $G=\Aut(X)$ and let $G_0 = \Aut_0(X)$ be the group of type-preserving
automorphisms of $X$.  \begin{enumerate}
\item\label{i:thomas} If there are $s,t \in S$ such that $q_s > 2$ and $m_{st} = \infty$ then
$G_0$ and $G$ are both nondiscrete.  \item\label{i:other} If all $q_s = 2$, or for all $t\in S$ with $q_{t} > 2$ we have $m_{st} = 2$ for
all $s \in S - \{t\}$, then $G_0$ is discrete, and $G$ is nondiscrete if and only if the nerve $L$ of
$(W,S)$ is flexible.
\end{enumerate} \end{theorem}

\noindent Case \eqref{i:thomas} of Theorem \ref{t:autom_gp_intro} follows easily from results of \cite{T1}, while if all $q_s = 2$ then the building $X$ is the Davis complex for $(W,S)$, in which case this result is due to Haglund--Paulin~\cite{HP1} and White~\cite{W}.

The second main result of Section~\ref{s:applications} is:

\begin{theorem}\label{t:transitive} Let $G$ be the automorphism group of a regular right-angled building $X$.  Then the
action of $G$ on $X$ is strongly transitive.
\end{theorem}

\noindent A group $G$ is said to act \emph{strongly transitively} on a building $X$ if it acts transitively on the set of
pairs $(\phi,\Sigma)$, where $\phi$ is a chamber of $X$, and $\Sigma$ is an apartment of $X$ containing $\phi$
(see Section~\ref{ss:rabs}).   By a theorem of Tits (see~\cite{D}), if $X$ is a thick building, meaning that there is some $q_s > 2$, then the group $G$ has a $BN$--pair.  For example, Bourdon's building $I_{p,q}$ is thick for all $q > 2$.  In fact, our techniques show that $\Stab_{G_0}(Y)$ acts transitively on all apartments containing a convex subcomplex $Y$ of $X$.  Theorem~\ref{t:transitive} was sketched for the case $X = I_{p,q}$ by Bourdon in~\cite[Proposition 2.3.3]{B}.  The general result will not surprise experts and, as pointed out by a referee, follows from, for example, results in \cite[Lemma 4.1]{DO} via methods as in the proofs of \cite[Theorems 4.2 and 4.4]{DO}.  

We would like to thank Indira Chatterji and Benson Farb for advice and encouragement, Kenneth S. Brown,
G. Christopher Hruska, Shahar Mozes, and Boris Okun for helpful conversations, Martin Bridson, Karen Vogtmann and an anonymous referee for
valuable comments on this manuscript, and the University of Chicago, MSRI and
Cornell University for supporting travel by both authors.

\section{Background}\label{s:background}

In Section~\ref{ss:lattices} we briefly describe the natural topology on $G$ the automorphism group of a
locally finite polyhedral complex $X$ and characterize uniform lattices in $G$.  We present some
necessary background on Coxeter groups and Davis complexes in Sections~\ref{ss:coxeter_groups}
and~\ref{ss:Davis_complexes} respectively, then discuss right-angled buildings in
Section~\ref{ss:rabs}.  Next in Section~\ref{ss:complexes_of_groups} we recall the basic theory of
complexes of groups and use this to construct the standard uniform lattice $\G_0$ in the automorphism
group of a regular right-angled building $X$.  Finally, Section~\ref{ss:covering_theory} contains
necessary definitions and results from covering theory for complexes of groups, and
Section~\ref{ss:group_actions} recalls the theory of group actions on complexes of groups.

\subsection{Lattices for polyhedral complexes}\label{ss:lattices}

Let $G$ be a locally compact topological group.  Recall that a discrete subgroup $\G \leq G$ is a \emph{lattice}
if $\G \bs G$ carries a finite $G$--invariant measure, and that $\G \leq G$ discrete is a \emph{uniform}
lattice if $\G \bs G$ is compact.

Let $X$ be a connected, locally finite polyhedral complex, and let $G=\Aut(X)$ be the group of automorphisms, or
cellular isometries, of $X$.  Then $G$, equipped with the compact-open topology, is a locally compact
topological group. In this topology, a countable neighborhood basis of the identity in $G$ consists of automorphisms
which fix larger and larger combinatorial balls in $X$ (we give the definition of combinatorial balls for $X$ a right-angled building in Section~\ref{sss:rabs} below).  A subgroup $\G$ of $G$ is discrete
if and only if, for each cell $\s$ of $X$, the stabilizer $\G_\s$ is a finite group.  Using a normalization of
the Haar measure on $G$ due to Serre~\cite{S}, and by the same arguments as for tree lattices (see Chapter 1 of~\cite{BL}),
if $G\bs X$ is compact, then $\G \leq G$ is a uniform lattice in $G$ exactly when $\G$ acts cocompactly
on $X$ with finite cell stabilizers.

\subsection{Coxeter groups}\label{ss:coxeter_groups}

We recall some necessary definitions and results.  Our notation and terminology in this section mostly follow
Davis~\cite{D}.

A \emph{Coxeter group} is a group $W$ with a finite generating set $S$ and presentation of the
form \[ W = \langle s\in S \mid (s t )^{m_{st}} = 1\rangle \] where $m_{ss} = 1$ for all $s \in S$, and if
$s \neq t$ then $m_{st}$ is an integer $\geq 2$ or $m_{st} = \infty$, meaning that there is no relation
between $s$ and $t$.   The pair $(W,S)$ is called a \emph{Coxeter system}.

Given a Coxeter system $(W,S)$, a \emph{word} in the generating set $S$ is a finite sequence \[\mathbf{s} =
(s_1,\ldots,s_k)\] where each $s_i \in S$.  We denote by $w(\mathbf{s})= s_1
\cdots s_k$ the corresponding element of $W$.  A word $\mathbf{s}$ is said to be \emph{reduced} if the element
$w(\mathbf{s})$ cannot be represented by any shorter word.  Tits proved that a word $\mathbf{s}$ is reduced
if and only if it cannot be shortened by a sequence of operations of either deleting a subword of the form $(s,s)$, or
replacing an alternating subword $(s,t,\ldots)$ of length $m_{st}$ by the alternating word $(t,s,\ldots)$ of the
same length $m_{st}$ (see Theorem 3.4.2~\cite{D}).  In particular, this implies:

\begin{lemma}\label{l:coxeter}  Any word in $S$ representing some $w\in W$ must involve all of the elements of $S$ that are used in any reduced word representing $w$. \end{lemma}

A Coxeter group $W$, or a Coxeter system $(W,S)$, is said to be
\emph{right-angled} if all $m_{st}$ with $s \neq t$ are equal to $2$ or $\infty$. That is, in a right-angled
Coxeter system, every pair of generators either commutes or has no relation.

\begin{examples}\label{egs:racg} Many later definitions and constructions will be illustrated by the following examples of
right-angled Coxeter groups.  \begin{enumerate}  \item\label{e:free_product}[free product] Let $W$ be the free
product of $n$ copies of $\Z / 2\Z$.  Then $W$ is a right-angled Coxeter group with presentation \[W =
\langle s_1,\ldots,s_n \mid s_i^2 = 1 \rangle.\]  In particular, if $n = 2$, then $W$ is the infinite
dihedral group. \item\label{e:mixed_dimension}[mixed dimension] Let $W$ be the free product of $\Z/ 2\Z$ with the direct
product $(\Z/2\Z \times \Z / 2\Z)$.  Then $W$ is a right-angled Coxeter group with presentation \[W =
\langle s_1, s_2, s_3 \mid s_i^2 = 1, (s_2 s_3)^2 = 1 \rangle.\] \item\label{e:hexagon}[hyperbolic hexagon] Let $W$ be the group generated by reflections in the sides of a regular right-angled hyperbolic hexagon.  Then $W$ is a right-angled Coxeter group with presentation \[ W = \langle s_1, \ldots,s_6 \mid s_i^2 = 1, (s_i s_{i+1})^2
= 1 \rangle\] where the subscripts of the $s_i$ are numbered cyclically. \end{enumerate}\end{examples}

\subsection{Davis complexes}\label{ss:Davis_complexes}

Let $(W,S)$ be a Coxeter system (not necessarily right-angled).  In this section we recall the construction
of the Davis complex $\Sigma$ for $(W,S)$, mostly following~\cite{D}.

For each subset $T$ of $S$, we define by $W_T:=\langle T \rangle$ the \emph{special subgroup} of $W$ generated by the elements $s\in T$.   By convention, $W_\emptyset$ is the trivial group.  A subset $T$ of $S$ is \emph{spherical} if $W_T$
is finite, in which case we say that $W_T$ is a \emph{spherical special subgroup}.  Denote by $\cS$ the set of
all spherical subsets of $S$.  Then $\cS$ is partially ordered by inclusion.  The poset $\cS_{> \emptyset}$
is an abstract simplicial complex, denoted by $L$ and called the \emph{nerve} of $(W,S)$.  In other words, the
vertex set of $L$ is $S$, and a nonempty set $T$ of vertices spans a simplex $\s_T$ in $L$ if and only if $T$ is
spherical.

\begin{examples} The nerves $L$ of Examples~\ref{egs:racg} above are as follows.
\begin{enumerate}
\item\label{dummy_label1}[free product] The $n$ vertices $\{ s_1 \}, \ldots, \{s_n\}$, with no higher-dimensional simplices.
\item\label{dummy_label2}[mixed dimension] A vertex $\{s_1\}$, and an edge joining the vertices $\{s_2\}$ and $\{s_3\}$.
\item\label{dummy_label3}[hyperbolic hexagon] A hexagon with vertices labeled cyclically $\{s_1\},\ldots,\{s_6\}$.
\end{enumerate}
\end{examples}

We denote by $K$ the geometric realization of the poset $\cS$.  Equivalently, $K$ is the cone on the barycentric
subdivision of the nerve $L$ of $(W,S)$.  Note that $K$ is compact and contractible, since it is the cone on a
finite simplicial complex.  Each vertex of $K$ has \emph{type} a spherical subset of $S$, with the cone point
having type $\emptyset$.

For each $s \in S$ let $K_s$ be the union of the (closed) simplices in $K$ which contain the vertex $\{s\}$ but
do not contain the cone point.  In other words, $K_s$ is the closed star of the vertex $\{s\}$ in the barycentric
subdivision of $L$.  Note that $K_s$ and $K_t$ intersect if and only if $m_{st}$ is finite. The family $(K_s)_{s
\in S}$ is a \emph{mirror structure} on $K$, meaning that $(K_s)_{s\in S}$ is a family of closed subspaces of
$K$, called \emph{mirrors}.  We call $K_s$ the \emph{$s$--mirror} of $K$.

\begin{lemma}[Lemma 7.2.5,~\cite{D}]\label{l:mirrors} Let $(W,S)$ be a Coxeter system and let $K$ be the
geometric realization of the poset $\cS$ of spherical subsets. \begin{enumerate} \item For each spherical
subset $T$, the intersection of mirrors $\cap_{s \in T} K_s$ is contractible. \item For each nonempty
spherical subset $T$, the union of mirrors $\cup_{s \in T} K_s$ is contractible. \end{enumerate} \end{lemma}

\noindent For any spherical subset $T$ of $S$, we call the intersection of mirrors $\cap_{s \in T} K_s$ a
\emph{face} of $K$, and the \emph{center} of this face is the unique vertex of $K$ of type $T$.  In
particular, the center of the $s$--mirror $K_s$ is the vertex $\{s\}$.

For each $x \in K$, put \[S(x):= \{ s
\in S \mid x \in K_s \}.\] Now define an equivalence relation $\sim$ on the set $W \times K$ by $(w,x) \sim
(w',x')$ if and only if $x = x'$ and $w^{-1}w' \in W_{S(x)}$.  The \emph{Davis complex} $\Sigma$ for $(W,S)$ is
then the quotient space: \[ \Sigma := (W \times K) / \sim. \] The types of vertices of $K$ induce types of
vertices of $\Sigma$, and the natural $W$--action on $W \times K$ descends to a type-preserving action on
$\Sigma$.

We identify $K$ with the subcomplex $(1,K)$ of $\Sigma$.  Then $K$, as well as any one of its translates by an
element of $W$, will be called a \emph{chamber} of $\Sigma$.  The subcomplexes $K_s$ of $K$, or any of their
translates by elements of $W$, will be called the \emph{mirrors} of $\Sigma$, and similarly for faces.

\begin{examples} For Examples~\ref{egs:racg} above:
\begin{enumerate} \item\label{dummy_label4}[free product] As shown in Figure~\ref{f:Davis_complex_tree} for $n =
3$, the chamber $K$ is the cone on $n$ vertices.  The Davis complex $\Sigma$ is
the barycentric subdivision of the $n$--regular tree, and its mirrors are the
midpoints of the edges of this tree.   If $n = 2$ then $\Sigma$ is homeomorphic
to the real line. \begin{figure}[ht] \scalebox{0.7}{\includegraphics{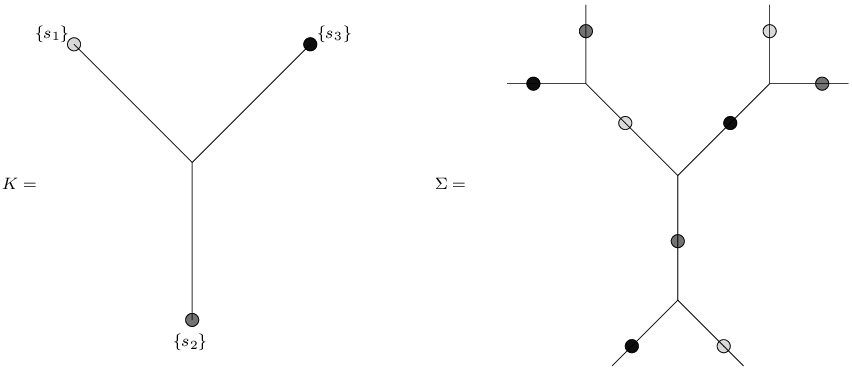}}
\caption{The chamber $K$ and the Davis complex $\Sigma$ for $W$ the free product
of $n = 3$ copies of $\Z/2\Z$.} \label{f:Davis_complex_tree} \end{figure} \item\label{dummy_label5}[mixed dimension]
The Davis complex $\Sigma$ for this example is sketched in Figure~1.2
of~\cite{D}; $\Sigma$ has both one-dimensional and two-dimensional maximal cells. \item\label{dummy_label6}[hyperbolic hexagon] The Davis complex $\Sigma$ for this example is homeomorphic
to the barycentric subdivision of the tesselation of the hyperbolic plane by
regular right-angled hexagons.  The mirrors are the edges of these hexagons.
\end{enumerate}
\end{examples}

\subsection{Right-angled buildings}\label{ss:rabs}

We first discuss general chamber systems and buildings in Section~\ref{sss:chambers}, before specializing to
the right-angled case in Section~\ref{sss:rabs}.  The local structure of right-angled buildings, which is
important for our proofs, is described in Section~\ref{sss:local}.  Again, we mostly follow Davis~\cite{D}.

\subsubsection{Chamber systems and buildings}\label{sss:chambers}

A \emph{chamber system} over a set $S$ is a set $\Phi$ of \emph{chambers} together with a family of
equivalence relations on $\Phi$ indexed by the elements of $S$.  For each $s \in S$, two chambers are
\emph{$s$--equivalent} if they are equivalent via the equivalence relation corresponding to $s$; they are
\emph{$s$--adjacent} if they are $s$--equivalent and not equal.  Two chambers are \emph{adjacent} if they are
$s$--adjacent for some $s \in S$.  A \emph{gallery} in $\Phi$ is a finite sequence of chambers
$(\phi_0,\ldots,\phi_k)$ such that $\phi_{j-1}$ is adjacent to $\phi_{j}$ for $1 \leq j \leq k$. A chamber
system is \emph{gallery-connected} if any two chambers can be connected by a gallery.  The \emph{type} of a
gallery $(\phi_0,\ldots,\phi_k)$ is the word $\mathbf{s} = (s_1,\ldots,s_k)$, where $\phi_{j-1}$ is
$s_j$--adjacent to $\phi_j$ for $1 \leq j \leq k$, and a gallery is \emph{minimal} if its type is a reduced
word.

\begin{definition} For $(W,S)$ a Coxeter system, the \emph{abstract Coxeter complex} $\mathbf{W}$ of $W$ is
 the chamber system with chambers the elements of $W$, and two chambers $w$ and $w'$ being $s$--adjacent, for
 $s \in S$, if and only if $w' = ws$.\end{definition}

\begin{definition}
Suppose that $(W,S)$ is a Coxeter system.  A \emph{building of type $(W,S)$} is a chamber system $\Phi$ over
$S$ such that:
\begin{enumerate}
\item for all $s \in S$, each $s$--equivalence class contains at least two chambers; and
\item there exists a \emph{$W$--valued distance function} $\delta:\Phi \times \Phi \to W$, that is, given a
reduced word $\mathbf{s}=(s_1,\ldots,s_k)$, chambers $\phi$ and $\phi'$ can be joined by a
gallery of type $\mathbf{s}$ in $\Phi$ if and only if $\delta(\phi,\phi') = w(\mathbf{s}) = s_1\cdots s_k$.
\end{enumerate}
\end{definition}

Let $\Phi$ be a building of type $(W,S)$.  Then $\Phi$ is \emph{spherical} if $W$ is finite.  The
building $\Phi$ is \emph{thick} if for all $s \in S$, each $s$--equivalence class of chambers contains at
least three elements; a building which is not thick is \emph{thin}.  The building $\Phi$ is \emph{regular}
if, for all $s \in S$, each $s$--equivalence class of chambers has the same number of elements.

\begin{example} The abstract Coxeter complex $\mathbf{W}$ of $W$ is a regular thin building, with $W$--distance function $\delta$
given by $\delta(w,w') = w^{-1}w'$.
\end{example}

Suppose $\Phi$ is a building of type $(W,S)$.  An \emph{apartment} of $\Phi$ is an image of the abstract
Coxeter complex $\mathbf{W}$, defined above, under a map $\mathbf{W} \to \Phi$ which preserves
$W$--distances.  The building $\Phi$ has a \emph{geometric realization}, which we denote by $X$, and by abuse
of notation we call $X$ a \emph{building of type $(W,S)$} as well.  By definition of the geometric
realization, for each chamber of $\Phi$, the corresponding subcomplex of $X$ is isomorphic to the chamber $K$
defined in Section~\ref{ss:Davis_complexes} above, and for each apartment of $\Phi$, the corresponding
subcomplex of the building $X$ is isomorphic to the Davis complex $\Sigma$ for $(W,S)$.  The copies of
$\Sigma$ in $X$ are referred to as the \emph{apartments} of $X$, and the copies of $K$ in $X$ are the
\emph{chambers} of $X$.  Note that each vertex of $X$ thus inherits a type $T$ a spherical subset of $S$. The
copies of $K_s$, $s \in S$, in $X$ are the \emph{mirrors} of $X$, so that two chambers in $X$ are
$s$--adjacent if and only if their intersection is a mirror of type $s$.  The \emph{faces} of $X$ are its
subcomplexes which are intersections of mirrors.  Each face has type $T$ a spherical subset of $S$, and a
face of type $T$ contains a unique vertex of type $T$, called its \emph{center}.

The building $X$ may be metrized as follows:

\begin{theorem}[Davis, Moussong, cf. Theorems 18.3.1 and 18.3.9 of \cite{D}]\label{t:metrize}  Let $(W,S)$ be
a Coxeter system and let $X$ be a building of type $(W,S)$. \begin{enumerate} \item \label{t:buildCAT0} The
building $X$ may be equipped with a piecewise Euclidean structure, such that $X$ is a complete $\CAT (0)$
space. \item \label{t:buildCAT-1} The building $X$ can be equipped with a piecewise hyperbolic structure which is $\CAT (-1)$ if
and only if $(W,S)$ satisfies Moussong's Hyperbolicity Condition: \begin{enumerate} \item  there is no subset
$T\subset S$ such that $W_T$ is a Euclidean reflection group of dimension $\geq 2$; and \item  there is no
subset $T \subset S$ such that $W_T = W_{T'} \times W_{T''}$ for nonspherical subsets $T', T'' \subset S$.
\end{enumerate} \end{enumerate} \end{theorem}

\noindent Unless stated otherwise, we equip buildings $X$ with the $\CAT(0)$ metric of Part \eqref{t:buildCAT0} of
Theorem~\ref{t:metrize}.

\subsubsection{Right-angled buildings}\label{sss:rabs}

In this section we specialize to right-angled buildings.  A building $X$ of type $(W,S)$ is
\emph{right-angled} if $(W,S)$ is a right-angled Coxeter system.   Note that part \eqref{t:buildCAT-1} of
Theorem~\ref{t:metrize} above implies that a piecewise hyperbolic $\CAT (-1)$ structure exists for a right-angled
building $X$ if and only if the nerve $L$ has no squares without diagonals (``satisfies the no--$\Box$
condition'').

The following result classifies regular right-angled buildings.

\begin{theorem}[Proposition 1.2, \cite{HP}]\label{t:rabeu} Let $(W,S)$ be a right-angled Coxeter system and
$\{q_s\}_{s \in S}$ a family of cardinalities.  Then, up to isometry, there exists a unique building $X$ of
type $(W,S)$, such that for all $s \in S$, each $s$--equivalence class of $X$ contains $q_s$ chambers.
\end{theorem}

\noindent In the $2$--dimensional case, this result is due to Bourdon~\cite{B}.  According to~\cite{HP},
Theorem~\ref{t:rabeu} was proved by M. Globus, and was known also to M. Davis, T. Janusz\-kiewicz, and J.
\Swiatkowski.  We will refer to a right-angled building $X$ as in Theorem~\ref{t:rabeu} as a \emph{building
of type $(W,S)$ and parameters $\{q_s\}$}.   In Section~\ref{ss:complexes_of_groups} below, we recall a
construction, appearing in Haglund--Paulin~\cite{HP}, of regular right-angled buildings $X$ as universal
covers of complexes of groups.

The following definition will be important for our proofs below.

\begin{definition} Let $X$ be a building of type $(W,S)$.  Fix $K$ some chamber of $X$.  We define the
\emph{combinatorial ball} $Y_n$ of radius $n$ in $X$ inductively as follows.  For $n = 0$, $Y_0 = K$, and for
$n \geq 1$, $Y_n$ is the union of $Y_{n-1}$ with the set of chambers of $X$
which have nonempty intersection with $Y_{n-1}$.\end{definition}

\begin{examples} \begin{enumerate} \item Let $(W,S)$ be the free product of $n$ copies of $\Z/2\Z$, as in
part~\eqref{e:free_product} of Examples~\ref{egs:racg} above. For $1 \leq i \leq n$ let $q_{i} = q_{s_i} \geq
2$ be a positive integer. Then the right-angled building $X$ of type $(W,S)$ and parameters $\{q_i\}$ is a
locally finite tree.  Each mirror $K_i = K_{s_i}$ is a vertex of $X$ of valence $q_i$.  The remaining
vertices of $X$ are the centers of chambers and have valence $n$.  If $n = 2$ then $X$ is the barycentric
subdivision of the $(q_1,q_2)$--biregular tree, and each chamber of $X$ is the barycentric subdivision of an
edge of this tree. Figure~\ref{f:rab_tree} depicts the combinatorial ball $Y_2$ of radius $2$ in $X$ for an
example with $n = 3$.   \begin{figure}[ht] \scalebox{1.0}{\includegraphics{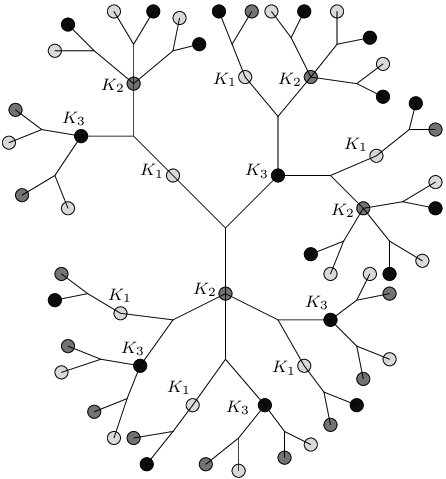}} \caption{The
combinatorial ball $Y_2$ of radius $2$, and mirrors contained in it, in the building $X$ of type $(W,S)$ and
parameters $q_1 = 2$, $q_2 = 4$ and $q_3 = 3$, where $W$ is the free product of $n = 3$ copies of $\Z/2\Z$.}
\label{f:rab_tree} \end{figure} \item In low dimensions, there are right-angled buildings $X$ which are also
hyperbolic buildings, meaning that their apartments are isometric to a (fixed) tesselation of hyperbolic
space $\Hyp^n$.  For this, let $P$ be a compact, convex, right-angled polyhedron in $\Hyp^n$; such polyhedra
$P$ exist only for $n \leq 4$, and this bound is sharp (Potyagailo--Vinberg~\cite{PV}).  Let $(W,S)$ be the
right-angled Coxeter system generated by reflections in the codimension one faces of $P$, and let $X$  be a
building of type $(W,S)$. By Theorem~\ref{t:metrize} above, $X$ may be equipped with a piecewise hyperbolic
structure which is $\CAT(-1)$.  Moreover, in this metric the apartments $\Sigma$ of $X$ are the barycentric
subdivision of the tesselation of $\Hyp^n$ by copies of $P$.  Thus $X$ is a hyperbolic building. For example,
Bourdon's building $I_{p,q}$ (see~\cite{B}) is of type $(W,S)$ and parameters $\{q_s\}$, where $W$ is
generated by reflections in the sides of $P$ a regular right-angled hyperbolic $p$--gon ($p \geq 5$), and each
$q_s = q \geq 2$.  Figure~\ref{f:oneballIpq} below shows the
combinatorial ball $Y_1$ of radius $1$ in $X = I_{6,3}$.\end{enumerate} \end{examples}

\begin{figure}[ht]
\scalebox{0.6}{\includegraphics{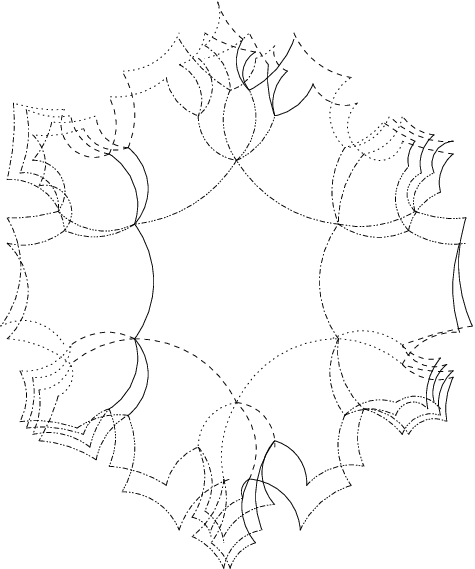}}
\caption{The combinatorial ball $Y_1$ of radius $1$ in Bourdon's building
$I_{6,3}$.}
\label{f:oneballIpq}
\end{figure}

\subsubsection{Local structure of right-angled buildings}\label{sss:local}

In our proofs below, we will rely on the following observations concerning the links of vertices in right-angled
buildings.

Let $X$ be a regular right-angled building of type $(W,S)$ and parameters $\{q_s\}_{s \in S}$.  Suppose $\s$
is a vertex of $X$, of type a \emph{maximal} spherical subset $T$ of $S$.  Then the link of $\s$ in $X$,
denoted by $\Lk_\s(X)$, is the (barycentric subdivision of the) join of $|T|$ sets of points, denoted
$V_t$, of cardinalities $|V_t| = q_t$ for each $t \in T$. For example, the link of each vertex of Bourdon's building
$I_{p,q}$ is the complete bipartite graph $K_{q,q}$, which may be thought of as the join of $2$ sets of $q$
points.  In fact, $\Lk_\s(X)$ is a (reducible) spherical building, of type $(W_T,T)$.

Now consider $\phi$ a chamber of $X$ such that the vertex $\s$ is in $\phi$.  Denote by $k_\phi$ the
subcomplex of the link $\Lk_\s(X)$ corresponding to simplices in $X$ which are contained in the chamber
$\phi$.  For example, in Bourdon's building $I_{p,q}$, the subcomplex $k_\phi$ is an edge of the graph $K_{q,q}$.  By abuse
of terminology, we call $k_\phi$ a \emph{maximal simplex} of $\Lk_\s(X)$.   (This is justified by recalling
that the chamber $\phi=K$ is the cone on the barycentric subdivision of the nerve $L$, hence $\phi$ is
homeomorphic to the cone on $L$.  Moreover, the maximal simplices of $L$ correspond precisely to the
maximal spherical subsets of $S$.)

Two chambers $\phi$ and $\phi'$ of $X$ containing $\s$ are adjacent in
$X$ if and only if the corresponding maximal simplices $k_\phi$ and $k_{\phi'}$ in $\Lk_\s(X)$ share a
codimension one face in $\Lk_\s(X)$.  Hence, a gallery of chambers in $X$, each chamber of which contains
$\s$, corresponds precisely to a gallery of maximal simplices in the spherical building $\Lk_\s(X)$.

\subsection{Basic theory of complexes of groups}\label{ss:complexes_of_groups}

In this section we sketch the theory of complexes of groups, due to Haefliger~\cite{BH}.  The sequence of
examples in this section constructs the regular right-angled building $X$ of Theorem~\ref{t:rabeu} above, as
well as the standard uniform lattice $\Gamma_0$ in $\Aut(X)$.  We postpone the definitions of
morphisms and coverings of complexes of groups to Section~\ref{ss:covering_theory} below.  All references
to~\cite{BH} in this section are to Chapter III.$\cC$.

In the literature, a complex of groups $G(Y)$ is constructed over a space or set $Y$ belonging to
various different categories, including simplicial complexes, polyhedral complexes, or, most
generally, \emph{scwols} (small categories without loops):

\begin{definition}\label{d:scwol} A \emph{scwol} $X$ is the disjoint
union of a set $V(X)$ of vertices and a set $E(X)$ of edges, with each edge $a$ oriented from its
initial vertex $i(a)$ to its terminal vertex $t(a)$, such that $i(a) \not = t(a)$.  A pair of
edges $(a,b)$ is \emph{composable} if $i(a)=t(b)$, in which case there is a third edge $ab$, called
the \emph{composition} of $a$ and $b$, such that $i(ab)=i(b)$, $t(ab)=t(a)$, and if $i(a) = t(b)$
and $i(b)=t(c)$ then $(ab)c = a(bc)$ (associativity).  \end{definition}

\noindent We will always assume scwols are \emph{connected} (see Section~1.1,~\cite{BH}).

\begin{definition}\label{d:action_on_scwol} An \emph{action of a
group $G$ on a scwol $X$} is a homomorphism from $G$ to the group of automorphisms of the scwol
(see Section 1.5 of~\cite{BH}) such that for all $a \in E(X)$ and all $g \in G$: \begin{enumerate}
\item $g.i(a) \not = t(a)$; and \item \label{c:no_inversions} if $g.i(a) = i(a)$ then $g.a = a$.
\end{enumerate} \end{definition}

Suppose $X$ is a right-angled building of type $(W,S)$, as defined in Section~\ref{ss:rabs} above. Recall
that each vertex $\sigma \in V(X)$ has a type $T \in \cS$.  The edges $E(X)$ are then naturally oriented by
inclusion of type.  That is, the edge $a$ joins a vertex $\sigma$ of type $T$ to a vertex $\sigma'$
of type $T'$, with $i(a)=\sigma$ and $t(a)=\sigma'$, if and only if $T \subsetneq T'$.  It is clear that the
sets $V(X)$ and $E(X)$ satisfy the properties of a scwol. Moreover, if $Y$ is a subcomplex of $X$, then the
sets $V(Y)$ and $E(Y)$ also satisfy Definition~\ref{d:scwol} above.  By abuse of notation, we identify $X$ and $Y$
with the associated scwols.  Note that a group of type-preserving automorphisms of $X$ acts according to Definition~\ref{d:action_on_scwol}, and that if $G=\Aut(X)$ is not type-preserving we may replace $X$ by a barycentric subdivision, with suitably oriented edges, on which $G$ does act according to Definition~\ref{d:action_on_scwol}.

We now define complexes of groups over scwols.

\begin{definition}
A \emph{complex of groups} $G(Y)=(G_\sigma, \psi_a, g_{a,b})$ over a scwol
$Y$ is given by: \begin{enumerate} \item a group $G_\sigma$ for each
$\sigma \in V(Y)$, called the \emph{local group} at $\sigma$;
\item a monomorphism $\psi_a: G_{i(a)}\rightarrow G_{t(a)}$ along the edge $a$ for each
$a \in E(Y)$; and
\item for each pair of composable edges, a twisting element $g_{a,b} \in
G_{t(a)}$, such that \[ \Ad(g_{a,b})\circ\psi_{ab} = \psi_a
\circ\psi_b
\] where $\Ad(g_{a,b})$ is conjugation by $g_{a,b}$ in $G_{t(a)}$,
and for each triple of composable edges $a,b,c$ the following
\emph{cocycle condition} holds:
\[\psi_a(g_{b,c})\,g_{a,bc} = g_{a,b}\,g_{ab,c}.\] \end{enumerate}
\end{definition}

\noindent A complex of groups is \emph{simple} if each $g_{a,b}$ is trivial.

Let $X$ be a regular right-angled building of type $(W,S)$ and parameters $\{q_s\}_{s \in S}$, where each
$q_s$ is an integer $q_s \geq 2$.  We construct $X$ and the standard uniform lattice $\G_0 < \Aut(X)$ using a
simple complex of groups $G_X(Y_0)$, which we now define.

\begin{definition}[Compare~\cite{HP}, p. 160]\label{d:GXY_0}  Let $K = Y_0$ be the cone on the barycentric
subdivision of the nerve $L$ of $(W,S)$ (see Section~\ref{ss:Davis_complexes} above). The simple complex of
groups $G_X(Y_0)$ over $Y_0$ is defined as follows.  For each $s \in S$ let $G_s$ be the cyclic group $\Z /
q_s \Z$.  The local group at the vertex of type $\emptyset$ of $Y_0$ is the trivial group.  The local group
at the vertex of type $T$ a nonempty spherical subset of $S$ is defined to be the direct product
\[G_T:=\prod_{s \in T} G_s.\]  All monomorphisms between local groups are natural inclusions, and all
$g_{a,b}$ are trivial. \end{definition}

\noindent Figures~\ref{f:complex_of_groups_tree} and~\ref{f:complex_of_groups_Ipq} below show this complex of
groups for the right-angled Coxeter systems in Examples~\ref{egs:racg} above.

\begin{figure}[ht] \scalebox{0.7}{\includegraphics{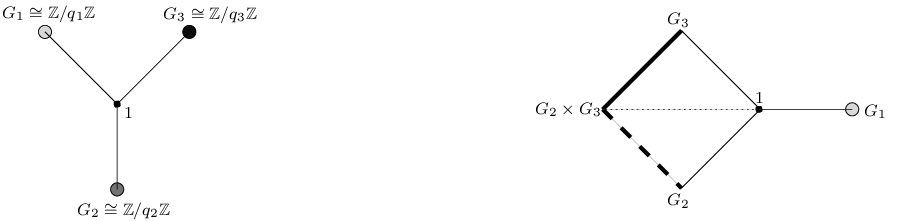}} \caption{The complex of groups
$G_X(Y_0)$ when $W$ is as in parts~\eqref{e:free_product} (on the left) and~\eqref{e:mixed_dimension} (on the right) of
Examples~\ref{egs:racg} above.  In both figures, $G_i = G_{q_i} $.} \label{f:complex_of_groups_tree} \end{figure}

\begin{figure}[ht]
\scalebox{1.0}{\includegraphics{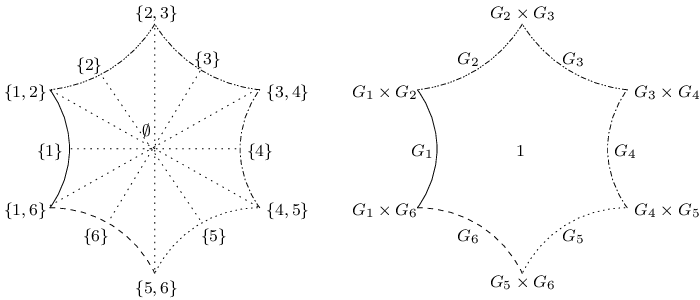}}
\caption{Types of vertices in $Y_0$, and the complex of groups $G_X(Y_0)$, for
Bourdon's building $I_{6,q}$.  Each group $G_i$ is isomorphic to $\Z/q\Z$.}
\label{f:complex_of_groups_Ipq}
\end{figure}

Suppose a group $G$ acts on a scwol $X$, as in Definition~\ref{d:action_on_scwol} above.  Then the
quotient $Y = G \bs X$ also has the structure of a scwol, and the action of $G$ on $X$ induces a complex of
groups $G(Y)$ over $Y$, as follows.  Let $p:X \to Y$ be the natural projection.  For each $\sigma \in V(Y)$, choose a lift $\overline\sigma \in V(X)$ with $p(\overline\sigma) = \sigma$.  The local group $G_\s$ of $G(Y)$ is then defined to be the the stabilizer of $\overline\sigma$ in $G$, and the monomorphisms $\psi_a$ and the elements $g_{a,b}$ are defined using further choices.  A complex of groups is \emph{developable} if it is isomorphic (see Definition~\ref{d:morphism} below) to a complex of groups $G(Y)$ induced by such an action.

Complexes of groups, unlike graphs of groups, are not in general developable.  We now discuss a sufficient condition for
developability.  Let $Y$ be a scwol equipped with the metric structure of a polyhedral complex.  An example is $Y$ a
subcomplex of a right-angled building $X$.    Each vertex $\sigma$ of $Y$ has a \emph{local development} in $G(Y)$,
which is, roughly speaking, a simplicial complex determined combinatorially by the cosets in $G_\s$ of the local groups
at vertices adjacent to $\s$.  The local group $G_\sigma$ acts naturally on the local development at $\s$, with quotient
the star of $\sigma$ in $Y$.  (The links of local developments for the complex of groups $G_X(Y_0)$ are described in the next example.) The metric on $Y$ induces a metric on the local development at $\s$.  We say that $G(Y)$
has \emph{nonpositive curvature} if, for every $\sigma \in V(Y)$, this induced metric on the local development at
$\sigma$ is locally $\CAT(0)$.  A nonpositively curved complex of groups $G(Y)$ is developable (Theorem
4.17,~\cite{BH}).

\begin{example}\label{e:developable} We continue the notation of Definition~\ref{d:GXY_0} above, and show that $G_X(Y_0)$ is
nonpositively curved and thus developable.  By Section 4.20 of~\cite{BH}, it is enough to check that the
local development at each vertex $\s$ of $Y_0$, of type $T$ a \emph{maximal} spherical subset of $S$, is
locally $\CAT(0)$.  By Gromov's Link Condition (see~\cite{BH}), for this, it suffices to show that the link
of the local development at $\s$ in $G_X(Y_0)$ is $\CAT(1)$. Now, for each proper subset $T'$ of $T$, there
is a unique vertex of $Y_0$ adjacent to $\s$ of type $T'$.  In particular, for each $t \in T$, there is a
unique vertex of $Y_0$ adjacent to $\s$ of type $T - \{t\}$.  It follows, by the construction of $G_X(Y_0)$
and Section 4.20 of~\cite{BH}, that the link of the local development at $\s$ is the join of $|T|$ sets of
points, of respective cardinalities $|G_T/G_{T - \{t\}}| = q_t$.  That is, the link of the local development
at $\s$ is the same as the link of a vertex of type $T$ in the building $X$.  As described in
Section~\ref{sss:local} above, the vertices of type $T$ in $X$ have links which are spherical buildings.
So these links are $\CAT(1)$.  Hence $G_X(Y_0)$ is nonpositively curved, and thus developable.\end{example}

The \emph{fundamental group $\pi_1(G(Y))$} of a complex of
groups $G(Y)$ is defined so that if $G(Y)$ is a simple complex of groups and $Y$ is simply
connected, then $\pi_1(G(Y))$ is isomorphic to the direct limit
of the family of groups $G_\sigma$ and monomorphisms $\psi_a$.

\begin{example} Since the chamber $Y_0=K$ is contractible, the fundamental group $\G_0 := \pi_1(G_X(Y_0))$ is
the graph product of the finite cyclic groups $(G_s)_{s \in S}$.  That is, $\G_0$ is the quotient of the free
product of the groups $(G_s)_{s \in S}$ by the normal subgroup generated by all commutators of the form
$[g_s,g_t]$ with $g_s \in G_s$, $g_t \in G_t$ and $m_{st} = 2$.\end{example}

If $G(Y)$ is a developable complex of groups, then it has a \emph{universal cover}
$\widetilde{G(Y)}$.  This is a connected, simply-connected scwol, equipped with an action of
$\pi_1(G(Y))$, so that the complex of groups induced by the action of the fundamental
group on the universal cover is isomorphic to $G(Y)$.  For each vertex $\s$ of $Y$, the star of any lift of $\s$ in $\widetilde{G(Y)}$ is isomorphic to the local development of $G(Y)$ at $\s$.

\begin{example}  By the discussion in Example~\ref{e:developable} above, the complex of groups $G_X(Y_0)$ is
developable.  By abuse of notation, denote by $X$ the universal cover of $G_X(Y_0)$.  Since the vertices of
$Y_0$ are equipped with types $T \in \cS$, the complex of groups $G_X(Y_0)$ is of \emph{type $(W,S)$} in the
sense defined in Section 1.5 of Gaboriau--Paulin~\cite{GP}.  As discussed above, the links of vertices of
$Y_0$ in their local developments are $\CAT(1)$ spherical buildings.  By an easy generalization of Theorem~2.1
of~\cite{GP}, it follows that the universal cover $X$ is a building of type $(W,S)$. (Section 3.3
of~\cite{GP} treats the case of right-angled hyperbolic buildings.)  By construction, the building $X$ is
regular, with each mirror of type $s$ contained in exactly $q_s = |G_s|$ distinct chambers.  Hence by
Theorem~\ref{t:rabeu} above, $X$ is the unique regular right-angled building of type $(W,S)$ and parameters
$\{ q_s \}$.
\end{example}

Let $G(Y)$ be a developable complex of groups over $Y$, with universal cover $X$ and fundamental
group $\G$.  We say that $G(Y)$ is \emph{faithful} if the action of $\G$ on $X$ is faithful, in which
case $\G$ may be identified with a subgroup of $\Aut(X)$.  If $X$ is locally finite, then with the
compact-open topology on $\Aut(X)$, by the discussion in Section~\ref{ss:lattices} above, the subgroup $\Gamma$ is discrete if and only if all local groups of $G(Y)$ are finite, and a discrete subgroup $\G$ is a uniform lattice in $\Aut(X)$ if and only if
the quotient $Y \cong \G \bs X$ is compact.

\begin{example} Since the local group in $G_X(Y_0)$ at the vertex of type $\emptyset$ of
$Y_0$ is
trivial, the fundamental group $\Gamma_0$ acts faithfully on the universal cover $X$.  Since $G_X(Y_0)$ is a complex of finite groups, $\Gamma_0$ is
discrete, and since $Y_0$ is compact, $\Gamma_0$ is a uniform lattice in $\Aut(X)$.
\end{example}

\noindent We call $\Gamma_0$ the \emph{standard uniform
lattice}.

\subsection{Covering theory for complexes of groups}\label{ss:covering_theory}

In this section we state necessary definitions and results from covering theory for complexes of groups.  As in Section~\ref{ss:complexes_of_groups} above, all
references to~\cite{BH} are to Chapter III.$\cC$.

We first recall the definitions of morphisms and coverings of complexes of groups.  In each of the
definitions below, $Y$ and $Z$ are scwols, $G(Y)=(G_\s,\psi_a)$ is a simple complex of groups over $Y$, and
$H(Z)=(H_\tau,\theta_a, h_{a,b})$ is a complex of groups over $Z$.  (We will only need morphisms and
coverings from simple complexes of groups $G(Y)$.)

\begin{definition}\label{d:morphism} Let $f: Y\to Z$ be a morphism of scwols (see Section~1.5
of~\cite{BH}).  A \emph{morphism} $\Phi: G(Y) \to H(Z)$ over $f$ consists of: \begin{enumerate}
\item a homomorphism $\phi_\sigma: G_\sigma \to H_{f(\sigma)}$ for each $\sigma \in V(Y)$, called
the \emph{local map} at $\s$; and
\item\label{i:commuting} an element $\phi(a) \in H_{t(f(a))}$ for each $a \in E(Y)$, such that the following diagram commutes
\[\xymatrix{
G_{i(a)}   \ar[d]^-{\phi_{i(a)}} \ar[rrr]^{\psi_a} & & & G_{t(a)} \ar[d]^-{\phi_{t(a)}}
\\
H_{f(i(a))}  \ar[rrr]^{\Ad(\phi(a))\circ \theta_{f(a)}} & & & H_{f(t(a))}
}\]
and for all pairs of
composable edges $(a,b)$ in $E(Y)$, \[ \phi(ab) = \phi(a) \,\psi_a(\phi(b))h_{f(a),f(b)} . \]
\end{enumerate} \end{definition}

\noindent A morphism is \emph{simple} if each element $\phi(a)$ is trivial.  If $f$ is an isomorphism of scwols, and each $\phi_\sigma$
an isomorphism of the local groups, then $\Phi$ is an \emph{isomorphism of complexes of groups}.

\begin{definition}\label{d:covering} A morphism $\Phi:G(Y) \to H(Z)$ over $f:Y\to Z$ is a
\emph{covering of complexes of groups} if further: \begin{enumerate}\item each $\phi_\sigma$ is
injective; and \item \label{i:covbijection} for each $\sigma \in V(Y)$ and $b \in E(Z)$ such that
$t(b) = f(\sigma)$, the map of cosets \[  \left(\coprod_{\substack{a \in f^{-1}(b)\\ t(a)=\sigma}} G_\sigma /
\psi_a(G_{i(a)})\right) \to H_{f(\sigma)} / \theta_b(H_{i(b)})\] induced by $g \mapsto
\phi_\sigma(g)\phi(a)$ is a bijection.\end{enumerate}\end{definition}

We will need the following general result on functoriality of
coverings, which is implicit in~\cite{BH}, and stated and proved explicitly in~\cite{LT}.

\begin{theorem}\label{t:coverings} Let $G(Y)$ and $H(Z)$ be complexes of groups over scwols $Y$ and $Z$ and let
$\Phi:G(Y) \to H(Z)$ be a covering of complexes of groups. If $G(Y)$ has nonpositive curvature
(hence is developable) then $H(Z)$ has nonpositive curvature, hence $H(Z)$ is developable.
Moreover, $\Phi$ induces a monomorphism of fundamental groups \[ \eta: \pi_1(G(Y)) \to
\pi_1(H(Z))\] and an $\eta$--equivariant isomorphism of universal covers \[
\widetilde{G(Y)} \to \widetilde{H(Z)}.\] \end{theorem}

See~\cite{LT} for the definition of an $n$--sheeted covering of complexes of groups, and the result that if $G(Y) \to H(Z)$ is an $n$--sheeted covering then the monomorphism $\eta:\pi_1(G(Y)) \to \pi_1(H(Z))$ in Theorem~\ref{t:coverings} above embeds $\pi_1(G(Y))$ as an index $n$ subgroup of $\pi_1(H(Z))$.

\subsection{Group actions on complexes of groups}\label{ss:group_actions}

The theory of group actions on complexes of groups was introduced in~\cite{T2}.  Let $G(Y)$ be a complex
of groups. An \emph{automorphism} of $G(Y)$ is an isomorphism $\Phi:G(Y) \to G(Y)$.  The set of all
automorphisms of $G(Y)$ forms a group under composition, denoted $\Aut(G(Y))$.  A group \emph{$H$ acts
on $G(Y)$} if there is a homomorphism $\rho:H \to \Aut(G(Y))$.  If $H$ acts on $G(Y)$, then in
particular $H$ acts on the scwol $Y$ in the sense of Definition~\ref{d:action_on_scwol} above, so we may
say that the $H$--action on $Y$ \emph{extends to an action on $G(Y)$}.  Denote by $\Phi^h$ the
automorphism of $G(Y)$ induced by $h \in H$.  We say that the \emph{$H$--action on $G(Y)$ is by simple
morphisms} if each $\Phi^h$ is a simple morphism.

\begin{theorem}[Thomas, Theorem 3.1 of~\cite{T2} and its proof]\label{t:group_action}  Let $G(Y)$ be a simple complex of groups over a connected scwol $Y$.  Suppose that the action of a group $H$ on $Y$ extends to an action by simple morphisms on $G(Y)$.
Then the $H$--action on $G(Y)$ induces a complex of groups $H(Z)$ over $Z = H \bs Y$, well-defined up to isomorphism of complexes of groups, such that:
\begin{itemize}
\item if $G(Y)$ is faithful and the $H$--action on $Y$ is faithful then $H(Z)$ is faithful;
\item there is a covering of complexes of groups $G(Y) \to H(Z)$; and
\item if $H(Z)$ is developable and $H$ fixes a point of $Y$, then $H \hookrightarrow \pi_1(H(Z))$.
\end{itemize}\end{theorem}

In particular, if the covering $G(Y) \to H(Z)$ is finite-sheeted, as occurs for example if $G(Y)$ is a complex of finite groups over a finite scwol $Y$, then $\pi_1(G(Y))$ is a finite index subgroup of $\pi_1(H(Z))$.

\section{Discreteness of normalizers}\label{s:norm}

Let $G$ be the group of automorphisms of a locally finite polyhedral complex $X$ (not necessarily a
building), and suppose $G \bs X$ is compact.  In this section we show that for any uniform lattice $\G$ of
$G$, the normalizer $N_G(\G)$ is discrete.  Recall from Section~\ref{ss:lattices} above that a uniform
lattice $\G$ in $G=\Aut(X)$ acts cocompactly on $X$, and fix a compact fundamental domain $D$ for this
action.

\begin{lemma}\label{l:ZDiscrete} The centralizer of $\G$ in $G$, denoted $Z_G(\G)$, is discrete in $G$.
\end{lemma} \begin{proof} Suppose otherwise.    Then there is a sequence $g_k \rightarrow \Id_X$ with $\Id_X
\neq g_k \in Z_G(\G)$.   Since $D$ is compact, it follows that for $k$ sufficiently large $g_k|_D=\Id_D$.
Let $x\in X$.  As $D$ is a $\G$--fundamental domain, $x\in \gamma D$ for some $\gamma \in \G$, that is,
$\gamma ^{-1} x \in D$.  It follows that $g_k(\gamma^{-1} x)=\gamma^{-1} x$, so $\gamma^{-1} g_k x=\gamma^{-1}
x$ as $g_k\in Z_G(\G)$.  Thus $g_kx=x$ for all $x\in X$ so $g_k=\Id_X$, a contradiction. \end{proof}

\begin{proposition}\label{p:Ndiscrete} The uniform lattice $\G$ is a finite index subgroup of its normalizer
$N_G(\G)$.  In particular, $N_G(\G)$ is discrete in $G$. \end{proposition} \begin{proof} By Lemma~\ref{l:ZDiscrete},
it follows directly from Proposition 6.2(c) of \cite{BL} that $N_G(\G)$ is also discrete.  Since $\G < N_G(\G)$, the
group $N_G(\G)$ is also a uniform lattice in $G$.  The ratio of covolumes of $\G$ and $N_G(\G)$ gives the index of
$\G$ in $N_G(\G)$.  In particular, this index is finite.  \end{proof}

We now sketch an alternative argument for $N_G(\G)$ being discrete, which was suggested to us by G.
Christopher Hruska, and which uses the theory of group actions on complexes of groups
(Section~\ref{ss:group_actions} above).  A uniform lattice $\G$ of $G=\Aut(X)$ is the fundamental group
of a complex of groups $G(Y)$, where $Y=\G \bs X$ is compact and the local groups of $G(Y)$ are finite.
Thus the group $\Aut(G(Y))$ of automorphisms of $G(Y)$ is a finite group.  Any element $g \in N_G(\G)$
induces an automorphism of $Y$, and this automorphism extends to an action on the complex of groups
$G(Y)$ (not necessarily by simple morphisms).   The induced action of $g$ on $G(Y)$ is trivial if and
only if $g \in \G$, so we have an isomorphism $N_G(\G) / \G \to \Aut(G(Y))$, hence $N_G(\G)$ is
discrete.

\section{Unfoldings}\label{s:unfolding}

We now introduce the technique of ``unfolding", which will be used in our proofs in Sections~\ref{s:Proof}
and~\ref{s:applications} below.  Let $X$ be a regular right-angled building. We first, in
Section~\ref{ss:clumps}, define \emph{clumps}, which are a class of subcomplexes of $X$ that includes the
combinatorial balls $Y_n \subset X$.  For each clump $\cC$ we then construct a canonical complex of groups
$G_X(\cC)$ over $\cC$, and we define a clump $\cC$ to be \emph{admissible} if $G_X(\cC)$ is developable with
universal cover $X$.  In Section~\ref{ss:unfoldingclumps}, we define the \emph{unfolding} of a clump $\cC$.
The main result of this section is Proposition~\ref{p:unfolding_admissible}, which shows that if $\cC$ is
admissible then any unfolding of $\cC$ is also admissible.  Finally, in Section~\ref{ss:unfoldingscover}, we
prove in Proposition~\ref{p:unfoldings_cover} that if $\cC$ is a clump obtained by a finite sequence of
unfoldings of the chamber $Y_0$, then there is a covering of complexes of groups $G_X(\cC) \to G_X(Y_0)$.  As
a corollary, we obtain a sequence $\G_n$ of uniform lattices in $G=\Aut(X)$, such that each $\G_n$ has
fundamental domain $Y_n$ and is of finite index in the standard uniform lattice $\G_0$.

\subsection{Complexes of groups over clumps}\label{ss:clumps}

Let $X$ be a regular right-angled building of type $(W,S)$.   In this section, we define clumps and, for each clump $\cC$ in $X$, construct a canonical complex of groups $G_X(\cC)$ over $\cC$.

We will say that two mirrors of $X$
are \emph{adjacent} if the face which is their intersection has type $T$ with $|T| = 2$. Since $(W,S)$ is
right-angled, it is immediate that:

\begin{lemma} If two adjacent mirrors are of the same type, then they are contained in adjacent chambers.  If
two adjacent mirrors are of different types, then there is a chamber of $X$ which contains both of these
mirrors.
\end{lemma}

\begin{definition}\label{d:clump} Let $X$ be a regular right-angled building of type $(W,S)$.\begin{itemize}\item A \emph{clump} in $X$ is a gallery-connected union of
chambers $\cC$ such that at least one mirror of $\cC$ is contained in only one chamber of $\cC$.  \item The
\emph{boundary} of a clump $\cC$, denoted $\partial \cC$, is the union of all the mirrors in $\cC$ each of
which is contained in only one chamber of $\cC$.  \item Two mirrors in $\partial\cC$ are \emph{type-connected} if
they are of the same type and are equivalent under the equivalence relation generated by adjacency.  \item If $\cC$
is a clump, then a maximal union of type-connected mirrors in $\partial \cC$ will be called a
\emph{type-connected component} or \emph{side} of $\cC$, and the \emph{type} of the side is the type of the
mirrors in the side. \end{itemize}\end{definition}

Let $\cC$ be a clump in $X$.  For a vertex $\s\in V(\cC)$ of type $T$, the \emph{boundary type} of $\s$ is
the subset \[ \{ s \in T \mid \mbox{an $s$--mirror containing $\s$ is contained in $\partial \cC$} \}.\] Note
that if $\s\notin \bC$, then the boundary type of $\s$ is $\emptyset$.

We now define a simple complex of
groups $G_X(\cC)$ over $\cC$.  For $s\in S$, let $G_s$ be the cyclic group $\Z/q_s\Z$.  For a vertex $\s$ in
$\cC$, we denote by $G_\s(\cC)$ the local group at $\s$ in $G_X(\cC)$.  Then $G_X(\cC)$ is defined as
follows: \begin{itemize} \item The local group $G_\s(\cC)$ at each vertex $\s \in \cC-\partial \cC$ is
trivial. \item The local group $G_\s(\cC)$ at a vertex $\s \in \partial \cC$ of boundary type $T$ is the
direct product \[ G_T:= \prod _{s\in T} G_s.\] \item The monomorphisms $\psi_a$ are natural
inclusions, for each edge $a$ in $\cC$. \item The twisting elements $g_{a,b}$ are all trivial. \end{itemize}

A clump $\cC$ is \emph{admissible} if $G_X(\cC)$ is developable and its universal cover is (isomorphic to)
$X$.  If $\cC$ is an admissible clump, then we may identify $\cC$ with a fundamental domain in $X$ for the
fundamental group of $G_X(\cC)$.  The \emph{preimage} or \emph{lift} of a vertex $\s \in V(X)$ in $\cC$ is the
unique vertex $\s'$ of $\cC$ which is in the same orbit as $\s$ under the action of the fundamental group of $G_X(\cC)$ on
$X$.  Lifts of edges and of chambers in $\cC$ are defined similarly.

\begin{example} The chamber $Y_0=K$ is an admissible clump since $G_X(Y_0)$ is precisely the defining complex
of groups for the standard uniform lattice $\G_0$ (see Section~\ref{ss:complexes_of_groups} above). \end{example}

\begin{example} Figure~\ref{f:nonadmissible} below depicts the complex of groups $G_X(\cC)$ over a clump
$\cC$ in the product $X=T_{q_s} \times T_{q_t}$ of regular trees of valences $q_s$ and $q_t$ respectively.
This clump is nonadmissible, since the link of the vertex $\s$ in the local development of $G_X(\cC)$ at $\s$
is not a complete bipartite graph, so this link is not the same as the link of a vertex in $X$. \end{example}

\begin{figure}[ht]
\scalebox{0.8}{\includegraphics{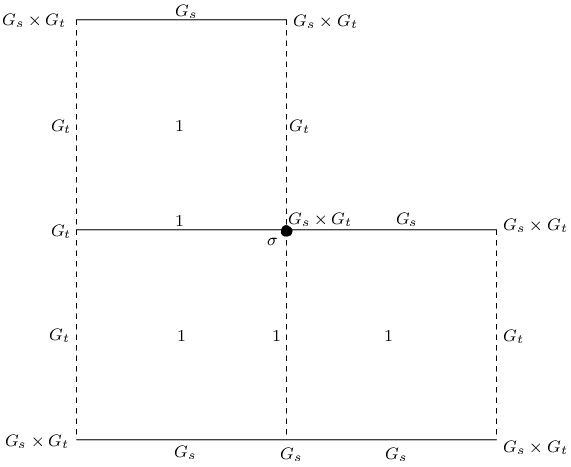}}
\caption{The complex of groups $G_X(\cC)$ with $\cC$ a nonadmissible clump.}
\label{f:nonadmissible}
\end{figure}

\subsection{Unfolding along a side of an admissible clump} \label{ss:unfoldingclumps}

Given an admissible clump, we now define a process, called \emph{unfolding}, that yields larger admissible
clumps.  In particular, as shown in Lemma~\ref{l:unfoldingYn} below, by starting with $Y_0$ and iterating
this process, one can obtain each of the combinatorial balls $Y_n$.  The main result of this section is
Proposition~\ref{p:unfolding_admissible} below, which shows that if $\cC$ is admissible then any unfolding of
$\cC$ is admissible.  Hence each $Y_n$ is admissible.

Let $\cC$ be an admissible clump in $X$, and let $\cK $ be a side of $\cC$ of type $u$.  The \emph{unfolding} of $\cC$
along $\cK $ is the clump \[U_{\cK}(\cC) := \cC \cup \{\mbox{chambers }\phi\subset X \mid \mbox{the }u\mbox{--mirror
of }\phi\mbox{ is contained in }\cK \}.\]

\begin{lemma}\label{l:unfoldingYn} The combinatorial ball $Y_n$ can be obtained by a sequence of unfoldings
beginning with a base chamber $Y_0$. \end{lemma}

\begin{proof} By induction, it suffices to show that the combinatorial ball $Y_n$ can be obtained from
$Y_{n-1}$ by a sequence of unfoldings.  Let $\cK_1$, $\cK_2, \dots, \cK_k$ denote the sides of $Y_{n-1}$.
First unfold $Y_{n-1}$ along $\cK_1$ to obtain a new clump $U_{\cK_1}(Y_{n-1})$.  For $i>1$, if $\cK_i$ does
not intersect $\cK_1$, then $\cK_i$ is also a side of $U_{\cK_1}(Y_{n-1})$.  Otherwise, replace $\cK_i$ by
the side of $U_{\cK_1}(Y_{n-1})$ containing $\cK_i$.  Then unfold along the (potentially extended) side
$\cK_2$.  Iterating this process, the clump $\cC$ obtained by unfolding along each of the (extended) sides
$\cK_1, \cK_2, \ldots, \cK_k$ is the combinatorial ball $Y_n$.  Figure~\ref{f:unfolding_Ipq} below
illustrates this process for obtaining $Y_1$ from $Y_0$ in Bourdon's building $I_{6,3}$. \end{proof}

\begin{figure}[ht]
\scalebox{0.8}{\includegraphics{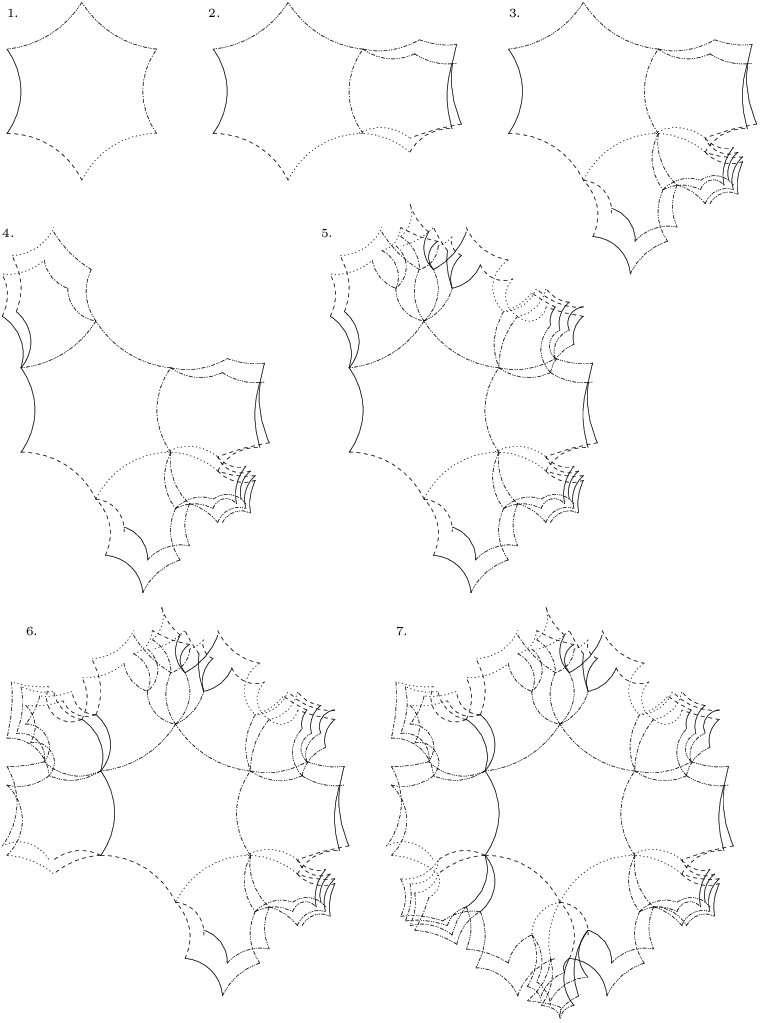}}
\caption{Unfolding $Y_0$ to get $Y_1$ in $I_{6,3}$.}
\label{f:unfolding_Ipq}
\end{figure}

We say that a vertex $\s$ in a clump $\cC \subset X$ is \emph{fully interior} if every chamber in $X$
containing $\s$ is in $\cC$.  Note that if some $q_s>2$, then a vertex can be in $\cC-\partial \cC$ without being
fully interior.  However, if $\cC$ is admissible, then since $X$ is the universal cover of $G_X(\cC)$, but
interior local groups in $G_X(\cC)$ are all trivial, it follows that every interior vertex of $\cC$ is fully
interior.

We call the local development at a vertex $\s$ in $G_X(\cC)$ \emph{complete} if it is the same as the local development
of a vertex of the same type in $G_X(Y_0)$, that is, if it is the star of $\s$ in $X$. We note that:

\begin{lemma}\label{l:complete_one_chamber} If $\s$ is a vertex of $\cC$ such that $\s$ is contained in only one
chamber of $\cC$, then the local development of $G_X(\cC)$ at $\s$ is complete.
\end{lemma}

We next prove several lemmas which will be used in this section and in Section~\ref{ss:unfoldingscover} below.

\begin{lemma} \label{l:boundary_type} Let $\cC$ be an admissible clump and let $\s\in \partial \cC$ be a
vertex of type $T$ and boundary type $T_{\bC}$.  If $s\in T_{\bC}$, then every mirror of type $s$ in $\cC$
containing $\s$ is actually contained in $\bC$, and there are exactly $\displaystyle \prod_{t\in T-T_{\bC}}
q_t$ such mirrors. \end{lemma}

\begin{proof} Since $\cC$ is admissible, the local development of $G_X(\cC)$ at $\s$ is complete, so the link
$\Lk_{\s}(\cC)$ of $\s$ in $\cC$ is the quotient of the link $\Lk_{\s}(X)$ of $\s$ in $X$ by the action of
the local group $G_{T_{\bC}}$.  Now, as discussed in Section~\ref{sss:local} above, $\Lk_\s(X)$ is the join
of $|T|$ sets of vertices $V_t$ for $t\in T$, of cardinalities respectively $|V_t| = q_t$.  By construction
of $G_X(\cC)$, the action of the local group $ G_{T_{\bC}}=\prod_{t\in T_{\bC}} G_t$  on $\Lk_\s(X)$ is
transitive on each set $V_t$ with $t\in T_{\bC}$, and is trivial on the sets $V_t$ for $t\notin T_{\bC}$.  It
follows that $\Lk_{\s}(\cC)$ is also a join of $|T|$ sets of vertices: it is the join of a singleton set for
each $t\in T_{\bC}$, along with the sets $V_t$ for $t\notin T_{\bC}$.  For each $s \in T_{\bC}$, the faces in
$\Lk_\s(\cC)$ corresponding to the $s$--mirrors of $\cC$ which contain $\s$ are precisely those faces in
$\Lk_\s(\cC)$ which are a join of $|T|-1$ vertices: the singleton sets corresponding to each $t\in
T_{\bC}-\{s\}$, together with one vertex from each of the sets $V_t$ for $t\notin T_{\bC}$.  There are
$\displaystyle \prod_{t\in T-T_{\bC}} q_t$ such faces.  Now, by construction of $G_X(\cC)$, a face $k_s$ in
$\Lk_{\s}(\cC)$ of type $s\in T_{\bC}$ corresponds to a mirror in the boundary of $\cC$ if and only if its
stabilizer in $G_{T_{\bC}}$ is nontrivial.  Since the action of $G_{T_{\bC}}$ fixes each vertex in the sets
$V_t$ for $t\notin T_{\bC}$, it follows that all such mirrors must be on the boundary of $\cC$. \end{proof}

Note that Lemma~\ref{l:boundary_type} implies that for an admissible clump $\cC$, the boundary type of a vertex $\s$
of type $T$ is actually equal to $\displaystyle\{ s \in T \mid \mbox{all $s$--mirrors containing $\s$ are contained in
$\partial \cC$} \}.$  This is not necessarily true in nonadmissible clumps.  For example, in
Figure~\ref{f:nonadmissible}, $s$ and $t$ are in the boundary type of $\s$ even though neither every $s$-- nor every
$t$--mirror in $\cC$ containing $\s$ is contained in $\bC$.

Suppose $\cC'$ is an admissible clump, that $\cK$ is a side of $\cC'$ of type $u$, and that
$\cC=U_\cK(\cC')$.  Let $\s$ be a vertex in $\partial \cC$ of type $T$ and let $\s'$ be the lift of $\s$ to
$\cC'$.

\begin{lemma}\label{l:PrecedingBoundaryType} If $T_{\partial\cC'}$ denotes the boundary type of $\s'$ in $\cC'$ and $T_{\partial\cC}$ is the boundary type of $\s$ in $\cC$, then $T_{\partial\cC'}\subset T_{\partial\cC}$.
\end{lemma}
\begin{proof}
Suppose $s\in T-T_{\partial\cC}$.  Then there are at least two $s$--adjacent chambers in $\cC$ whose intersection contains $\s$.  The lifts of these chambers to $\cC'$ are then $s$--adjacent chambers in $\cC'$ containing $\s'$.  Hence $s\in T - T'_{\partial\cC'}$.
\end{proof}

Let $\Ch_\cK$ denote the set of chambers in $\cC=U_\cK(\cC')$ that are not also in $\cC'$, that is, $\Ch_\cK$ is
the set of ``new chambers'' in $\cC$.  A \emph{sheet} of chambers in $\Ch_\cK$ is an equivalence class of chambers
under the equivalence relation generated by $S-\{u\}$ adjacency in $\Ch_\cK$.  So two chambers in $\Ch_\cK$ are in the
same sheet if and only if there is a gallery of chambers in $\Ch_\cK$ such that the type of each adjacency is in
$S-\{u\}$.

\begin{lemma}\label{l:sheets} If $\cK$ is a side of $\cC'$ of type $u$, there are $q_u-1$ sheets in $\Ch_\cK$.
\end{lemma}

\begin{proof} Choose $K_u \subset \cK$ a mirror of type $u$.  There are $q_u-1$ chambers in $\Ch_\cK$
glued along $K_u$.  Call these chambers $\phi_1, \phi_2, \ldots, \phi_{q_u-1}$.  Since $\cK$ is type-connected, any
$\phi \in \Ch_\cK$ is in the same sheet as some $\phi_i$.  Now suppose there are $1 \leq i \neq j \leq q_u - 1$ such
that $\phi_i$ and $\phi_j$ are in the same sheet.  Then there is a gallery of chambers in $\Ch_\cK$ from $\phi_i$ to
$\phi_j$ with the type of each consecutive adjacency being an element of $S-\{u\}$.  The chambers $\phi_i$ and $\phi_j$ are
$u$--adjacent, since they are both glued to the mirror $K_u$, so the sequence ($\phi_i$, $\phi_j$) is also a gallery
in $X$.  By the definition of a building, and, more specifically, using the $W$--valued distance function, it follows
that $u$ is equal to a product of elements in $S-\{u\}$.  This is a contradiction, since $u \notin W_{S-\{u\}}$.
Hence there are exactly $q_u-1$ sheets in $\Ch_\cK$, namely the equivalence classes of each of $\phi_1$, $\phi_2$,
\ldots, $\phi_{q_u-1}$. \end{proof}

We now prove the main result of this section, that unfolding preserves admissibility.

\begin{proposition} \label{p:unfolding_admissible} Let $\cC_0$ be an admissible clump in $X$.  If $\cC$ is a
clump obtained from $\cC_0$ through a finite sequence of unfoldings, then $\cC$ is an admissible clump.   \end{proposition}

\begin{proof} By induction, it suffices to show that if $\cC'$ is an admissible clump and $\cK$ is a side of
$\cC'$ of type $u$, then the clump \[\cC= U_{\cK}(\cC')\] is admissible, that is, that $G_X(\cC)$ is
developable with universal cover $X$.  We will show that for each maximal spherical subset $T\subset S$, the local
development at each vertex $\s \in \cC$ of type $T$ is complete.  It will then follow that $G_X(\cC)$ is
developable with universal cover $X$, by similar arguments to those used for $G_X(Y_0)$ in
Section~\ref{ss:complexes_of_groups} above.

Let $\s$ be a vertex of $\cC$ of type $T$ a maximal spherical subset of $S$.  If $\s \in \cC' - \cK$, then
the set of chambers in $\cC$ containing $\s$ is the same as the set of chambers in $\cC'$ containing $\s$.
Thus the local development of $G_X(\cC)$ at $\s$ is the same as that of $G_X(\cC')$ at $\s$, since the
neighboring local groups are also all the same in the two complexes of groups.   Hence by induction the local
development at $\s$ is complete.

Thus it remains to consider the local developments of vertices in the side $\cK$ of $\cC'$ and in
$\cC-\cC'$.  We consider separately the three cases: \begin{description} \item[Case 1] $\s \in \cC - \cC'$
\item[Case 2] $\s \in \cK - K \cap \partial \cC$ \item[Case 3] $\s \in \cK \cap \partial \cC$
\end{description} as depicted in Figure~\ref{f:unfolding_cases} below.

\begin{figure}[ht] \scalebox{0.8}{\includegraphics{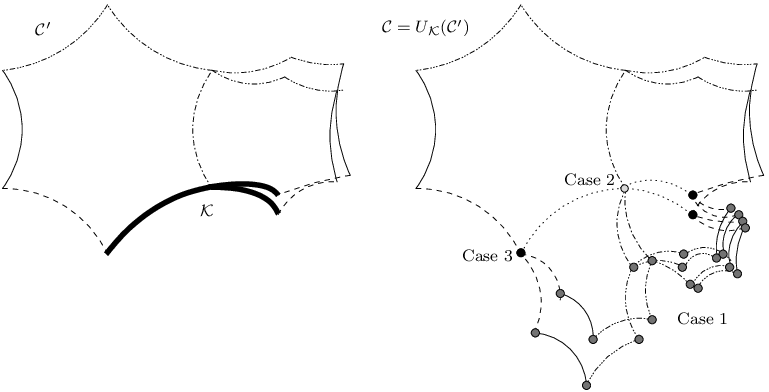}} \caption{A clump $\cC'$ in $I_{6,3}$, a side
$\cK$ of $\cC'$, the unfolding $\cC=U_{\cK}(\cC')$, and the three cases for vertices on the boundary of $\cC$ used in
the proof of Proposition~\ref{p:unfolding_admissible}.} \label{f:unfolding_cases} \end{figure}

\begin{description}

\item[Case 1]  Suppose $\s \in \cC - \cC'$ is a vertex of type $T$, and let $T_{\bC}$ be the boundary type of
$\s$ in $\cC$.  If $\s$ is contained in only one chamber of $\cC$, then by Lemma~\ref{l:complete_one_chamber}
above we are done.  Otherwise, we first prove:


\begin{lemma}\label{l:center} Let $\s \in \cC - \cC'$ and suppose $\s$ is contained in more than one chamber
of $\cC$.  Then there is a unique vertex of type $T - T_{\bC}$ adjacent to $\s$ in $\cC$. \end{lemma}

\begin{proof} Let $\phi_1$ and $\phi_2$ be two chambers of $\cC$ which contain $\s$.  Let $\phi_1'$ and $\phi_2'$ be
the lifts of $\phi_1$ and $\phi_2$ respectively to $\cC'$, and let $\s'$ be the lift of $\s$ to $\cC'$.

Since $\cC'$ is admissible, the link $\Lk_{\s'}(\cC')$ of $\s'$ in $\cC'$ is a join. Therefore there is a gallery
$\beta'$ in $\cC'$ from $\phi_1'$ to $\phi_2'$, of type say $\mathbf{t'}$, such that each chamber of the gallery
$\beta'$ contains the vertex $\s'$.  Without loss of generality, we may assume that $\beta'$ is a minimal gallery.
Thus $\mathbf{t'}$ is a reduced word.  The vertex $\s'$ is of type $T$, and every chamber in the gallery
$\beta'$ contains $\s'$, hence every letter in the reduced word $\mathbf{t'}$ must be an element of $T$.

The two chambers $\phi_1'$ and $\phi_2'$ both contain a mirror in the side $\cK$.  Since $\cK$ is a type-connected
component of mirrors of $\cC'$, there is also a gallery $\alpha'$ in $\cC'$ from $\phi_1'$ to $\phi_2'$, such that every
chamber in the gallery $\alpha'$ contains a mirror in the side $\cK$.  Let $\mathbf{s'}$ be the type of the
gallery $\alpha'$.  Since $\cK$ is of type $u$, it follows that every letter in $\mathbf{s'}$ commutes with $u$.

We now have two galleries $\alpha'$ and $\beta'$ from $\phi_1'$ to $\phi_2'$ in $\cC'$, of respective types
$\mathbf{s'}$ and $\mathbf{t'}$.  By Lemma~\ref{l:coxeter} above, since $\mathbf{t'}$ is reduced, every letter in
$\mathbf{t'}$ appears in $\mathbf{s'}$.  Hence every letter in $\mathbf{t'}$ is contained in $T$ and commutes with
$u$.

Now every letter  in the type $\mathbf{t'}$ of $\beta'$ commutes with $u$,  and the initial chamber
$\phi_1$ of $\beta'$ contains a mirror in the side $\cK$.  So by induction, every chamber in the gallery
$\beta'$ contains a mirror in the side $\cK$.

We claim that every letter in $\mathbf{t'}$ is actually contained in $T - T_{\bC}$, where $T_{\bC}$ is the boundary
type of $\s$ in $\cC$.  So suppose there is some $t \in T_{\bC}$ such that $t$ appears in the reduced word
$\mathbf{t'}$.  Denote by $T_{\bC'}$ the boundary type of $\s'$ in $\cC'$.  By Lemma~\ref{l:PrecedingBoundaryType},
we have that $T_{\bC'} \subset T_{\bC}$.  So assume first that $t \in T_{\bC'}$.  By Lemma~\ref{l:boundary_type}, since $\cC'$ is
admissible, every mirror of $\cC'$ of type $t$ which contains $\s'$ is in the boundary $\bC'$.  But the gallery
$\beta'$ is contained in $\cC'$, and every chamber in $\beta'$ contains the vertex $\s'$, so the gallery $\beta'$
cannot cross any mirror of type $t$ which also contains $\s'$.  So $t$ cannot be contained in $T_{\bC'}$.

We now have $t \in T_{\bC} - T_{\bC'}$.  Since $t \in T_{\bC}$, by definition there must be some chamber
$\widetilde\phi$ of $\cC$ which contains $\s$, such that the $t$--mirror of $\widetilde\phi$ is only
contained in one chamber of $\cC$.  Let $\overline\phi$ be a chamber of $X$ which is $t$--adjacent to
$\widetilde\phi$, and note that $\overline\phi$ is not in $\cC$.  Let $\widetilde\phi'$ be the lift of
$\widetilde\phi$ to $\cC'$.  Since $\widetilde\phi$ is in $\Ch_\cK$, the chambers $\widetilde\phi$ and
$\widetilde\phi'$ are $u$--adjacent.  Since $\cC'$ is admissible and $t \notin T_{\bC'}$, there is a chamber
say $\hat\phi'$ of $\cC'$ such that $\hat\phi'$ is $t$--adjacent to $\widetilde\phi'$.  Now, the letter $t$
commutes with $u$, and $\widetilde\phi'$ has its $u$--mirror contained in the side $\cK$ of $\cC'$.  Hence
the chamber $\hat\phi'$ of $\cC'$ also has its $u$--mirror contained in the side $\cK$.  Consider the gallery
$(\overline\phi,\widetilde\phi, \widetilde\phi', \hat\phi')$ in $X$.  This gallery has type $(t,u,t)$.
Since $t$ commutes with $u$, we have $tut = t^2 u = u$.  Hence $\overline\phi$ and $\hat\phi'$ are $u$--adjacent.
Therefore the $u$--mirror of $\overline\phi$ is contained in $\cK$.  But this implies that $\overline\phi$ is
in $\cC$, a contradiction.  We conclude that $t \in T - T_{\bC}$, as claimed.

We now have a minimal gallery $\beta'$ of type $\mathbf{t'}$ from $\phi_1'$ to $\phi_2'$ in $\cC'$, such that every chamber in the gallery $\beta'$ contains $\s'$, every chamber in $\beta'$ contains a mirror in the side $\cK$, and every letter in $\mathbf{t'}$ is contained in $T - T_{\bC}$ and commutes with $u$.

Next consider the gallery $\alpha$ from $\phi_1$ to $\phi_2$ obtained by concatenating the galleries
$(\phi_1,\phi_1')$, $\beta'$ and $(\phi_2',\phi_2)$.  Let $\mathbf{s}$ be the type of $\alpha$.  Then since
every letter in $\mathbf{t'}$ commutes with $u$, \[ w(\mathbf{s}) = u\, w(\mathbf{t'}) u = u^2 w(\mathbf{t'}) =
w(\mathbf{t'}).\] Since $\mathbf{t'}$ is a reduced word, it follows that there is a gallery, say $\beta$, in
$X$ from $\phi_1$ to $\phi_2$ of type $\mathbf{t'}$.  But every letter in $\mathbf{t'}$ commutes with $u$, so
every chamber in $\beta$ has a mirror contained in the side $\cK$.  Thus the gallery $\beta$ is contained in
$\cC$. That is, there is a minimal gallery $\beta$ from $\phi_1$ to $\phi_2$ in $\cC$, of type $\mathbf{t'}$,
such that every letter in $\mathbf{t'}$ is in $T - T_{\bC}$.

Let $\s_{T - T_{\bC}}^1$ and $\s_{T - T_{\bC}}^2$ be the vertices of types $T - T_{\bC}$ in $\phi_1$ and $\phi_2$
respectively.  Then since every letter in $\mathbf{t'}$ is in $T - T_{\bC}$, every chamber in the gallery
$\beta$ contains $\s_{T - T_{\bC}}^1$.  In particular, the chamber $\phi_2$ contains $\s_{T - T_{\bC}}^1$.
Hence $\s_{T - T_{\bC}}^1 = \s_{T - T_{\bC}}^2$.  We conclude that there is a unique vertex of type $T - T_{\bC}$ adjacent to $\s$ in $\cC$.  \end{proof}

By Lemma~\ref{l:center} and Lemma~\ref{l:boundary_type} above, the link $\Lk_\s(\cC)$ is a join of $|T-T_{\bC}|$ sets
of vertices $V_t$ of cardinality $q_t$ for each $t\in T-T_{\bC}$, and a singleton $\{v_s\}$ for each $s\in T_{\bC}$.
This is precisely the quotient of the link $\Lk_{\s}(X)$ of $\s$ in $X$ by the group $G_{T_{\bC}}$.  It follows that
the local development of $G_X(\cC)$ at $\s$ is complete.

\item[Case 2] Suppose $\s \in \cK - (\cK \cap \partial \cC)$.  Recall that the side $\cK$ has type $u$.  Let
$s\in S$ be in the boundary type of $\s$ in $\cC'$.  Then there is a mirror $K_s \subset \partial \cC'$ of
type $s$ containing $\s$.  If $s\neq u$, then $K_s \subset \partial \cC$, so $\s \in \partial \cC$, a
contradiction.  Hence the boundary type of $\s$ in $\cC'$ is $\{u\}$.  So the local group $G_\s(\cC')$ at
$\s$ in $G_X(\cC')$ is $G_u$.  Note that if the type of $\s$ in $X$ is also $\{u\}$, then $\s$ is the center
of a $u$--mirror in $\cK$, so all the chambers in $X$ containing $\s$ are in $\cC$, by definition of the
unfolding across $\cK$. Suppose then that the type of $\s$ is not $\{u\}$.  Let $\s_u$ be a vertex of type
$u$ in $\cC'$ that is adjacent to $\s$.  Since $\s_u$ is in $\cK$, the local group at $\s_u$ in $G_X(\cC')$
is also $G_u$, so in particular has index $1$ in the local group $G_\s(\cC')=G_u$.  By induction, the local
development at $\s$ in $G_X(\cC')$ is complete, so it follows that every vertex of type $u$ adjacent to $\s$
is in $\cC'$.  That is, every mirror of type $u$ containing $\s$ is in $\cC'$.  Thus every chamber of $X$
containing $\s$ is either in $\cC'$ or is adjacent to $\cC'$ along $\cK$.  Hence every such chamber is
contained in $\cC$, so $\s$ is fully interior in $\cC$, and it follows that the local development of
$G_X(\cC)$ at $\s$ is complete.

\item[Case 3] Suppose finally that $\s \in \cK \cap \partial \cC$ and let $T_{\bC'}$ be the boundary type of
$\s$ in $\cC'$.  Then the boundary type of $\s$ in $\cC$ is \hbox{$T_{\bC}=T_{\bC'}-\{u\}$} so its local
group in $G_X(\cC)$ is $G_{T_{\bC}}=G_{T_{\bC'}}/G_u$.  Now, since interior vertices of $\cC$ have trivial
local groups in $G_X(\cC)$, the number of chambers in the local development of
$G_X(\cC)$ at $\s$ is \[|G_{T_{\bC}}|\cdot \#\{ \mbox{chambers in }\cC\mbox{ containing }\s\}.\]  By
Lemma~\ref{l:boundary_type}, the number of chambers in the admissible clump $\cC'$ containing $\s$ is
$|G_{T-T_{\bC'}}|$. So by unfolding, we see that there are precisely $q_u\cdot |G_{T-T_{\bC'}}|$ chambers in
$\cC$ containing $\s$.  It follows that the number of chambers in the local development of $G_X(\cC)$ at $\s$
is precisely $|G_T|$.  In fact, we can describe the local structure at $\s$.

Since $\cC'$ is admissible, the link $\Lk_\s(\cC')$ of $\s$ in $\cC'$ is $G_{T_{\bC'}}\bs \Lk_{\s}(X)$.  This is the join
of the sets $V_t$ for $t\in T-T_{\bC'}$ and singletons $\{v_t\}$ for $t\in T_{\bC'}$.   Now
the local construction of $\cC$ from $\cC'$ at $\s$ consists of adding $q_u-1$ chambers along each $u$--mirror in
$\cK$ containing $\s$, so the link $\Lk_{\s}(\cC)$ of $\s$ in $\cC$ is as in Lemma~\ref{l:boundary_type}
above; it is the join of the $|T|$ sets of vertices $V_t$ for $t\in T-T_{\bC}$ and $\{v_t\}$ for $t\in T_{\bC}$.  It follows that the local development of $G_X(\cC)$ at $\s$ is
complete, as required.
\end{description} This completes the proof of Proposition~\ref{p:unfolding_admissible}. \end{proof}

\subsection{Unfoldings of $G_X(Y_0)$ cover $G_X(Y_0)$} \label{ss:unfoldingscover}

Recall from Section~\ref{ss:complexes_of_groups} above that the standard uniform lattice $\G_0$ is the
fundamental group of the complex of groups $G_X(Y_0)$ over a single chamber $Y_0$.  In this section, we show
that uniform lattices obtained via a sequence of unfoldings starting with $G_X(Y_0)$ are finite index
subgroups of $\G_0$.  The main result is the following proposition:

\begin{proposition}\label{p:unfoldings_cover}
Let $\cC_0 = Y_0$, and suppose that, for all $r > 0$, $\cC_r$ is a clump obtained by unfolding $\cC_{r-1}$ along a side $\cK_{r-1}$.  Then there is a covering of complexes of groups $G_X(\cC_r)\rightarrow G_X(\cC_0)$.  In particular, the fundamental group of $G_X(\cC_r)$ is a finite index subgroup of $\G_0$.
\end{proposition}

\noindent By Lemma \ref{l:unfoldingYn} above, the combinatorial balls $Y_n \subset X$ can be obtained by a
sequence of unfoldings of $Y_0$.  Let $\G_n$ be the fundamental group of $G_X(Y_n)$.  Then $\G_n$ is a
uniform lattice in $\Aut(X)$, and Proposition~\ref{p:unfoldings_cover} immediately implies:

\begin{corollary}\label{c:Gamma_n}
The lattices $\G_n$ are finite index subgroups of $\G_0$.
\end{corollary}

A key step in the proof of Proposition~\ref{p:unfoldings_cover} is provided by Proposition~\ref{p:labeling} below, the proof of which is at the end of this section.  It will be convenient to think of all groups $G_T$ for $T\subset S$ as natural subgroups of the direct product $\displaystyle G_S:=\prod_{s\in S} G_s$.

\begin{proposition}\label{p:labeling}
Let $\cC_r$ be as in Proposition \ref{p:unfoldings_cover} above.  Let  $p:\cC_r\rightarrow \cC_0$ be the natural morphism of scwols which sends a vertex of $\cC_r$ to the unique vertex of $\cC_0$ of the same type.  Then there is an edge labeling \[\lambda: E(\cC_r)\rightarrow G_S=\prod_{s\in S}G_s\] satisfying all of the following:
\begin{enumerate}
\item\label{i:one} $\lambda(a) \in G_{p(t(a))}$ for each $a\in E(\cC_r)$.
\item\label{i:two} For each pair of composable edges $(a,b)$ in $E(\cC_r)$,
\[ \lambda(ab) = \lambda(a)\lambda(b). \]
\item\label{i:three} For each $\s \in V(\cC_r)$ and $b \in E(\cC_0)$ such that $t(b) = p(\s)$, the map
\[\left(\coprod_{\substack{a \in p^{-1}(b)\\ t(a)=\s}} G_\s(\cC_r)
/ G_{i(a)}(\cC_r)\right) \to G_{p(\s)}(\cC_0) / G_{i(b)}(\cC_0)\] induced by $g
\mapsto g\lambda(a)$ is a bijection.
\end{enumerate}
\end{proposition}

\begin{proof}[Proof of Proposition \ref{p:unfoldings_cover}] We construct a covering $\Lambda:G_X(\cC_r) \to
G_X(\cC_0)$ over the natural morphism $p:\cC_r \to \cC_0$.  The local maps $\lambda_\s$ are defined to be the
identity map (if $\s$ is of type the empty set, or if the boundary type of $\s$ equals its regular type), or
natural inclusions (if $\s$ is an interior vertex of type $T$ not the empty set, or if the boundary type of
$\s$ is a proper subset of its regular type).  Note that the maps $\lambda_\s$ so defined are injective; by
abuse of notation, we write $g$ for $\lambda_\s(g)$.

We now use the edge labeling $\lambda$ provided by Proposition~\ref{p:labeling} above to complete the definition of $\Lambda$.  Since all local groups are abelian and the local maps $\lambda_\s$
are the identity or natural inclusions, the morphism diagram (see (2) of Definition \ref{d:morphism}) commutes no matter what the value of the $\lambda(a)$.  From the properties of $\lambda$ guaranteed by Proposition~\ref{p:labeling}, it thus follows that $\Lambda$ is a covering of complexes of groups. \end{proof}

\begin{proof}[Proof of Proposition \ref{p:labeling}] We proceed by induction on $r$ and write
$\lambda^r$ for the labeling of the edges of $\cC_r$.  See Figure~\ref{f:labeling} for an example. Given an edge $a\in E(\cC_r)$ such that $t(a)$ is a vertex of type $T$, we will choose an element $\lambda^r(a)$ of $G_T\subset G_S$.  Recall
that $G_T$ is the direct product of the cyclic groups $G_t$ for $t\in T$.  So, we can think of an
element of $G_T$ as an ordered $|T|$--tuple of elements of the cyclic groups $G_t$.  To define
$\lambda^r(a)$, it thus suffices to define elements $\lambda_t^r(a)\in G_t$ for each $t\in T$.  We
will refer to $\lambda_t^r(a)$ as the \emph{$t$--component} of $\lambda^r(a)$.

\begin{figure}[ht]
\scalebox{0.8}{\includegraphics{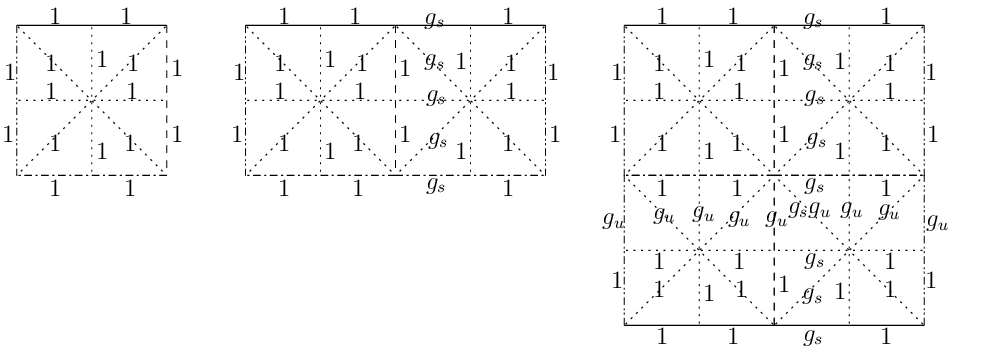}}
\caption{The labeling of edges in $\cC_0$, $\cC_1$ and $\cC_2$, where $\cC_0 = Y_0$ is the
barycentric subdivision of a square. Here $\cC_1$ is obtained from $\cC_0$ by
unfolding along a side of type $s$, with $G_s = \{ 1,g_s\}$, and $\cC_2$ is obtained from
$\cC_1$ by unfolding along a side of type $u$, with $G_u = \{1,g_u\}$.}
\label{f:labeling}
\end{figure}

To begin the induction, let $a\in E(\cC_0)$, with $t(a)$ of type $T$.  We define
$\lambda^0(a)$ to be the identity element in $G_T$.  Properties (1)--(3) in the statement of
Proposition~\ref{p:labeling} then hold trivially for $r=0$ with this labeling.

Suppose now that we inductively have a labeling $\lambda^{r-1}$ of $E(\cC_{r-1})$ satisfying
properties~\eqref{i:one}--\eqref{i:three} in the statement of Proposition~\ref{p:labeling}, and suppose $\cC_r$ is
obtained from $\cC_{r-1}$ by unfolding along a side $\cK=\cK_{r-1}$ of type $u$.  We first use the
labeling $\lambda^{r-1}$ to label the edges of $\cC_{r-1} \subset \cC_r$.  That is, for all $a\in
E(\cC_{r-1})$, define $\lambda^r(a):= \lambda^{r-1}(a)$.

Next, since $\cC_{r-1}$ is admissible, as in the proof of admissibility of unfoldings (Proposition
\ref{p:unfolding_admissible} above), we may think of $\cC_r$ as a subcomplex of the universal cover of
$G_X(\cC_{r-1})$.  For each edge $a\in E(\cC_r)-E(\cC_{r-1})$ there is then a unique preimage $a'\in
E(\cC_{r-1})$.  Define $\widehat\lambda^r(a):= \lambda^{r-1}(a')$.  If $t(a)\notin \cK$ then this is the
labeling we choose for $a$, that is, we set $\lambda^r(a):=\widehat\lambda^r(a)$.  If $t(a) \in \cK$ but
$i(a)\notin \cK$, then for $s\in S-\{u\}$ we define the $s$--component of $\lambda^r(a)$  to be the same as
that of $\widehat \lambda^r(a)$.  The $u$--component $\lambda_u^r(a)$ is then defined as follows.

Choose $K_u \subset \cK$ a mirror of type $u$ and let $c\in E(\cC_{r-1})$ be the edge with initial
vertex of type $\emptyset$ and terminal vertex the center of $K_u$.  Put $g=\lambda_u^{r-1}(c) \in
G_u$, that is, $g$ is the $u$--component of $\lambda^{r-1}(c)$.  There are $q_u-1$ chambers in
$\Ch_\cK$ glued along $K_u$.  Call these chambers $\phi_1, \phi_2, \ldots, \phi_{q_u-1}$.  For
each of these new chambers, we assign distinct elements of $G_u-\{g\}$, say $g_i$ is assigned to
$\phi_i$ for $1\leq i\leq q_u-1$, so that $G_u-\{g\}=\{g_i\mid 1\leq i\leq q_u-1\}$.  Now, for all
edges $a\in E(\phi_j)$ such that $t(a)\in K_u$ but $i(a)\notin K_u$, we define the $u$--component
$\lambda^r_u (a):=g_j$.  We then extend these $u$--components along the $q_u-1$ sheets of new
chambers described in Lemma \ref{l:sheets} above.  That is, for a chamber $\phi \in \Ch_\cK$ in
the same sheet as $\phi_j$, and for $a\in E(\phi)$ such that $t(a)\in \cK$ but $i(a)\notin \cK$,
we define $\lambda_u^r(a):=g_j$.

We must verify that this is well-defined.  Suppose $a\in E(\phi)\cap E(\phi')$ for some other $ \phi' \in
\Ch_\cK$.  We will show that $\phi$ and $\phi'$ are in the same sheet.  Consider the link $\Lk_{i(a)}(\cC_r)$
of $i(a)$ in $\cC_r$.  As in the proof of Proposition \ref{p:unfolding_admissible} above, since $\cC_r$ is
admissible, this is the join of sets of vertices.  In particular, the chambers $\phi$ and $\phi'$ correspond
to maximal simplices $k_\phi$ and $k_{\phi'}$ in this join.  A gallery in $\cC_r$ from $\phi$ to $\phi'$ and
containing $i(a)$ then corresponds to a sequence of maximal simplices in $\Lk_{i(a)}(\cC_r)$ from $k_\phi$ to
$k_{\phi'}$, which sequentially intersect along codimension one faces, that is, to a gallery in
$\Lk_{i(a)}(\cC_r)$. Such a sequence exists since $\Lk_{i(a)}(\cC_r)$ is a join.  Hence there is a gallery in
$\cC_r$ from $\phi$ to $\phi'$ each chamber of which contains the vertex $i(a)$.  Since $i(a) \notin \cK$,
this gallery cannot cross $\cK$.  It follows that $\phi$ and $\phi'$ are in the same sheet, as required.
Thus our assignment of the $u$--component of $\lambda^r(a)$, for edges $a \in \Ch_\cK$ with $t(a) \in \cK$
but $i(a) \notin \cK$, is well-defined.

This completes the definition of the labeling $\lambda^r$.  We now verify that $\lambda^r$
satisfies properties~\eqref{i:one}--\eqref{i:three} in the statement of Proposition~\ref{p:labeling}.

For~\eqref{i:one}, suppose $a\in E(\cC_r)$.  That $\lambda^r(a) \in G_{p(t(a))}$ follows immediately from the
above construction.

For~\eqref{i:two}, for each pair of composable edges $(a,b)$ in $E(\cC_r)$ we must show that $\lambda^r(ab) =
\lambda^r(a)\lambda^r(b)$.  If both $a$ and $b$ are in $E(\cC_{r-1})$, then this follows by induction.  Since
pairs of composable edges occur in chambers, the only other possibility is that $a$ and $b$ are edges in the
same chamber in $\Ch_\cK$.  It suffices to check that $\lambda^r_s(ab)=\lambda^r_s(a)\lambda^r_s(b)$ for all
$s\in S$. Let $a'$ and $b'$ be the preimages of $a$ and $b$ in $E(\cC_{r-1})$.  By induction,
$\lambda^{r-1}(a'b') = \lambda^{r-1}(a')\lambda^{r-1}(b')$.  The only possible difference between the labels
$\lambda^r(a)$ and $\lambda^{r-1}(a')$ is in the $u$--component, and similarly for $b$ and $ab$ (recall that
the side $\cK$ along which we unfolded is of type $u$).  Hence it suffices to show that
$\lambda^r_u(ab)=\lambda^r_u(a)\lambda^r_u(b)$.

By construction of $\lambda^r$, the only edges whose labels have different $u$--components from those of
their preimages are edges with terminal but not initial vertex in $\cK$.  For these edges, we have shown that
the $u$--component is determined by the chamber containing the edge.  Moreover, for a pair of composable
edges $(a,b)$, either none of $a$, $b$, and $ab$ have terminal but not initial vertex in $\cK$, or $ab$ and
exactly one of $a$ and $b$ do. In the latter case, by construction, the $u$--component of $ab$ is equal to
the $u$--component of the other edge ($a$ or $b$ but not both) with terminal but not initial vertex in
$\cK$.  It follows that $\lambda^r(ab) = \lambda^r(a)\lambda^r(b)$, as required.

Finally, for property~\eqref{i:three} in the statement of Proposition~\ref{p:labeling}, we show that for each
$\s \in V(\cC_r)$ and $b \in E(\cC_0)$ such that $t(b) = p(\s)$, the map of cosets \[\nu_r=\nu_{r,\s}:
\left(\coprod_{\substack{a \in p^{-1}(b)\\ t(a)=\s}} G_\s(\cC_r) / G_{i(a)}(\cC_r)\right) \to
G_{p(\s)}(\cC_0) / G_{i(b)}(\cC_0)\] induced by $g \mapsto g\lambda^r(a)$ is a bijection.  For this, we
assume that: \[ \mbox{$\s$ has type $T$ and $i(b)$ has type $U$}.\] So the codomain of $\nu_r$ is $G_T/G_U$,
and if $a \in p^{-1}(b)$ then $i(a)$ has type $U$.

Suppose $\s \in \cC_{r-1}-\cK$.  Then by induction and the construction of $\lambda^r$, $\nu_{r}=\nu_{r,\s}$
is bijective.

Suppose next that $\s \in \cK$.  Recall that $\cK$ is a side of type $u$.  If $T_{\partial\cC_{r-1}}$ is the
    boundary type of $\s$ in $\cC_{r-1}$ and $T_{\partial\cC_r}$ the boundary type of $\s$ in $\cC_r$,
    then $T_{\partial\cC_r} \cup \{ u \} = T_{\partial\cC_{r-1}}$.  Hence for all $\s \in \cK$,
    \begin{equation}\label{e:Gsigma} G_\s(\cC_r) \times G_u = G_\s(\cC_{r-1}).\end{equation} Denote by
    $p_{r-1}:\cC_{r-1} \to \cC_0$ and $p_r:\cC_r \to \cC_0$ the natural type-preserving morphisms of
    scwols.

    Assume first that $u \in U$.  Then by Lemma~\ref{l:boundary_type} above, since $\cC_{r-1}$ is
    admissible, \[ \{ a \in p_{r-1}^{-1}(b) \mid t(a) =\s \} = \{ a \in p_{r}^{-1}(b) \mid t(a) =\s \}
    \subset \cK.\] For all edges $a$ in this set, by construction $\lambda^r(a) = \lambda^{r-1}(a)$ and
    \begin{equation}\label{e:Gia}  G_{i(a)}(\cC_r) \times G_u = G_{i(a)}(\cC_{r-1}). \end{equation} By
    induction, the map $\nu_{r-1,\s}$ is bijective.  Therefore by Equations~\eqref{e:Gsigma}
    and~\eqref{e:Gia} it follows that $\nu_{r,\s}$ is bijective, as required.

    Now assume that $u \notin U$.  Then for all $a \in p_r^{-1}(b)$ with $t(a) = \s \in \cK$, we have
    $i(a) \notin \cK$.  Consider an edge $a' \in p_{r-1}^{-1}(b) \subset \cC_{r-1}$ with $t(a') = \s$.
    Since $i(a') \in \cC_{r-1} - \cK$, we now have  $G_{i(a')}(\cC_{r}) = G_{i(a')}(\cC_{r-1})$.  After
    unfolding, there are $q_u - 1$ images of $a'$ in $\Ch_\cK$, which we denote by $a_2,\ldots,a_{q_u}$.
    Put $a' = a_1$.  Then by construction, $G_u = \{\lambda^r_u(a_1), \lambda^r_u(a_2),\ldots,
    \lambda^r_u(a_{q_u})\}$, and for each $1 \leq j \leq q_u$ we have $G_{i(a_j)}(\cC_r) =
    G_{i(a')}(\cC_r)$.  Using Equation~\eqref{e:Gsigma} above, there is thus a natural bijection \[
    \zeta_{a'}: \left( \coprod_{j=1}^{q_u} G_\s(\cC_r)/G_{i(a_j)}(\cC_r) \right) \to
    G_\s(\cC_{r-1})/G_{i(a')}(\cC_{r-1})  \] induced by $g \mapsto g\lambda^r_u(a_j)$.  Note also that, by
    construction of the labeling $\lambda^r$, we have \begin{equation}\label{e:factors} \lambda^r(a_j) =
    \lambda^r_u(a_j) \lambda^{r-1}(a') g_j \end{equation} for some element $g_j \in G_u$.

    Let $\zeta$ be the disjoint union of the maps $\{ \zeta_{a'} \mid a' \in p_{r-1}^{-1}(b), t(a') = \s \}$.
Then $\zeta$ is a bijection from the domain of $\nu_{r,\s}$ to the domain of $\nu_{r-1,\s}$.  By induction
$\nu_{r-1,\s}$ is bijective. By Equation~\eqref{e:factors} above, the map $\nu_{r,\s}$ factors through
$\zeta$.  Hence $\nu_{r,\s}$ is bijective, as required.

    We have now proved that $\nu_{r,\s}$ is a bijection for all $\s \in \cC_{r-1}$.  For $\s \in
\cC_r-\cC_{r-1}$, let $\s'$ denote the unique preimage of $\s$ in $\cC_{r-1}$.

If $\s\in \cC_r-(\cC_{r-1} \cup \partial\cC_r)$, then the local structure at $\s$ in $\cC_r$ (meaning the set
of edges with terminal vertex $\s$, the local groups at the initial vertices of these edges, and the labels
of these edges) is identical to that at $\s'$ in $\cC_{r-1}$.   It follows by induction that $\nu_{r,\s}$ is
bijective.

    It remains to prove that $\nu_{r,\s}$ is bijective for $\s\in \partial \cC_r - (\partial\cC_r \cap
\cC_{r-1})$.  (Note that $\s$ is the same kind of vertex as in Case 1 in the proof of
Proposition~\ref{p:unfolding_admissible} above.)  Let $T_{\partial\cC_r}$ be the boundary type of $\s$ in
$\cC_r$.

    \begin{lemma} \label{l:U'well-defined}
     Suppose $a\in E(\cC_r)$ with $t(a) = \s$, and that $i(a)$ is of type $U'$ where
     \[U \subset U' \subset T.\]
     Then the boundary type $U'_{\partial\cC_r}$  of $i(a)$ in $\cC_r$ is given by \[U'_{\partial\cC_r}=U'\cap T_{\partial\cC_r}.\]  In particular, for all such edges $a$, the local group at $i(a)$ in $G_X(\cC_r)$ is the same.
    \end{lemma}

    \begin{proof} If $U'$ is the empty set then the boundary type $U'_{\partial\cC_r} \subset U'$ is also
    empty and we are done.  So suppose there is some $s\in U'$.  Then there is an $s$--mirror $K_s$ in
    $\cC_r$ which contains $i(a)$.  Since $U' \subset T$ and $t(a) = \s$, the mirror $K_s$ also contains
    $\s$.  By Lemma~\ref{l:boundary_type}, since $\cC_r$ is admissible, $s$ is in $ T_{\partial\cC_r}$ if and only if $K_s \subset
    \partial \cC_r$.  It follows that $s$ is in the boundary type of $i(a)$ if and only if $s$ is also in the
    boundary type of $\s$. \end{proof}

    \begin{lemma} \label{l:Unique_tau}
    Let \[U' = (T - T_{\partial\cC_r}) \cup U.\]
    Then there is a unique vertex of type $U'$ in $\cC_r$ adjacent to $\s$.
    \end{lemma}

    \begin{proof} Since $U'\subset T$ and $\cC_r$ is a gallery-connected union of chambers, there is at least
    one such vertex, say $\tau$.  By definition, the local group at $\s$ in $G_X(\cC_r)$ is
    $G_{T_{\partial\cC_r}}$ and the local group at $\tau$ in $G_X(\cC_r)$ is $G_{U'_{\partial\cC_r}}$.  Since
    $G_X(\cC_r)$ is admissible, there are thus \[
    \left|G_{T_{\partial\cC_r}}/G_{U'_{\partial\cC_r}}\right|=\prod_{s\in
    T_{\partial\cC_r}-U'_{\partial\cC_r}}q_s\] vertices of type $U'$ adjacent to $\s$ in $X$ that lift to
    $\tau$ in $\cC_r$.  But by admissibility of $G_X(\cC_0)$, the total number of vertices of type $U'$
    adjacent to $\s$ in $X$ is $\displaystyle \left|G_T/G_{U'}\right| = \prod_{s\in T-U'}q_s$.  Since by
    Lemma~\ref{l:U'well-defined} above \[T - U' =  T_{\partial\cC_r}-U'_{\partial\cC_r}\] it follows that
    $\tau$ is unique. \end{proof}

For a subset $R\subset S$, the \emph{projection} of an element $g\in G_S$ to $R$, or the \emph{$R$--projection} of $g$, is the projection of the ordered $|S|$--tuple $g$ to the components corresponding to $R$.  To simplify notation, write $p = p_r:\cC_r \to \cC_0$.

\begin{lemma}\label{l:coset_bijectivity_condition}
The map $\nu_{r,\s}$ is bijective if and only if the set of labels
\[ \{ \lambda^r(a) \mid a \in p^{-1}(b), t(a) = \s \} \]
has pairwise distinct projections to $T-(T_{\partial\cC_r} \cup U)$.
\end{lemma}

\begin{proof} By admissibility of $\cC_r$, the two sets \[ \left(\coprod_{\substack{a \in p^{-1}(b)\\
 t(a)=\s}} G_\s(\cC_r) / G_{i(a)}(\cC_r)\right)\quad\mbox{ and }\quad G_{p(\s)}(\cC_0) / G_{i(b)}(\cC_0)=G_T/G_U\] are
 finite sets of the same size.  So $\nu_r=\nu_{r,\s}$ is a bijection if and only if it is injective.  Note
 that $G_\s(\cC_r) = G_{T_{\partial\cC_r}}$ and that by Lemma~\ref{l:U'well-defined} above, for all $a \in
 p^{-1}(b)$ with $t(a) = \s$, we have $G_{i(a)}(\cC_r) = G_{U_{\partial\cC_r}} = G_{U \cap
 T_{\partial\cC_r}}$.

Let $a_1$ and $a_2$ be distinct edges in $p^{-1}(b)$ with $t(a_1) = t(a_2) = \s$.  Suppose $g, g' \in G_\s(\cC_r)$.  Now $\nu_r(gG_{i(a_1)}(\cC_r))=\nu_r(g'G_{i(a_2)}(\cC_r))$ if and only if $g\lambda^r(a_1)G_U=g'\lambda^r(a_2)G_U$.  Since $G_T$ is abelian, this equality of cosets holds if and only if $(g^{-1}g')\lambda^r(a_1)^{-1}\lambda^r(a_2) \in G_U$.

So if $\nu_r$ is injective, then in particular, putting $g'=1$, it follows that for all $g\in G_{T_{\partial\cC_r}}$, we have $\lambda^r(a_1)^{-1}\lambda^r(a_2) \notin gG_U$.  Hence $\lambda^r(a_1)^{-1}\lambda^r(a_2) \notin G_{T_{\partial\cC_r}\cup U}$.  That is, $\lambda^r(a_1)$ and $\lambda^r(a_2)$ have distinct $T-(T_{\partial\cC_r} \cup U)$ projections.

Conversely, suppose $\nu_r$ is not injective.  Then there are edges $a_1$ and $a_2$ and elements $g, g'\in
G_{T_{\partial\cC_r}}$ such that $\lambda^r(a_1)^{-1}\lambda^r(a_2)\in gg'^{-1}G_U$.  Then
$\lambda^r(a_1)^{-1}\lambda^r(a_2) \in G_{T_{\partial\cC_r}\cup U}$ and so the two labels $\lambda^r(a_1)$
and $\lambda^r(a_2)$ have the same
$T-(T_{\partial\cC_r} \cup U)$ projections. \end{proof}

Thus to prove that $\nu_r$ is a bijection, it suffices by Lemma~\ref{l:coset_bijectivity_condition} to show
that: for each pair of distinct edges $a_1, a_2 \in E(\cC_r)$ with $p(a_1)=p(a_2)=b$ and $t(a_1)=t(a_2)=\s$,
the labels $\lambda^r(a_1)$ and $\lambda^r(a_2)$ have distinct $T-(T_{\partial\cC_r} \cup U)$ projections.
Let $a_1'$, $a_2'$, and $\s'$ be the lifts of $a_1$, $a_2$, and $\s$, respectively, to $\cC_{r-1}$, and let
$T'_{\partial\cC_{r-1}}$ be the boundary type of $\s'$ in $\cC_{r-1}$.  By induction and Lemma
\ref{l:coset_bijectivity_condition}, the $T-(T'_{\partial\cC_{r-1}} \cup U)$ projections of
$\lambda^{r-1}(a_1')$ and $ \lambda^{r-1}(a_2')$ are distinct.  Since $\s\notin \cK$, the labels
$\lambda^r(a_1)$ and $\lambda^r(a_2)$ are the same as the labels $\lambda^{r-1}(a_1')$ and
$\lambda^{r-1}(a_2')$, respectively.  Hence the $T-(T'_{\partial\cC_{r-1}} \cup U)$ projections of
$\lambda^r(a_1)$ and $\lambda^r(a_2)$ are distinct.

 Now let $\tau$ be the unique vertex of type $U'=(T-T_{\partial\cC_r})\cup U$ in $\cC_r$ adjacent to $\s$, as
guaranteed by Lemma \ref{l:Unique_tau} above.  Let $d$ be the edge of $\cC_r$ with $i(d)=\tau$ and
$t(d)=\s$.  Since $U\subset U' \subset T$, there are edges $c_1$ and $c_2$ of $\cC_r$ such that
$i(c_1)=i(a_1)$ and $i(c_2)=i(a_2)$ are vertices of type $U$, and $t(c_1)=t(c_2)=\tau$ is of type $U'$.  We then have
compositions of edges $a_1=dc_1$ and $a_2=dc_2$, so by the already proved property (2) of the labeling
$\lambda^r$, we find that $\lambda^r(a_1)=\lambda^r(d)\lambda^r(c_1)$ and
$\lambda^r(a_2)=\lambda^r(d)\lambda^r(c_2)$.  Thus $\lambda^r(a_1)\lambda^r(a_2)^{-1} =
\lambda^r(c_1)\lambda^r(c_2)^{-1} \in G_{U'}$.  Note that by definition of $U'$ and
Lemma~\ref{l:U'well-defined}, $T_{\partial\cC_r} \cap U' = T_{\partial\cC_r} \cap U = U_{\partial\cC_r}$.  So
$\lambda^r(a_1)$ and $\lambda^r(a_2)$ have the same $T - U' = T_{\partial\cC_r}-U_{\partial\cC_r}$
projections.  Since they have different $T-(T'_{\partial\cC_{r-1}}\cup U)$ projections, it follows from Lemma
\ref{l:PrecedingBoundaryType} above that they have different $T-(T_{\partial\cC_r} \cup U)$ projections, as
required.

This completes the proof of Proposition~\ref{p:labeling}.\end{proof}

\section{Proof of the Density Theorem}\label{s:Proof}

We are now ready to complete the proof of the Density Theorem.  The main results we use are those of
Sections~\ref{ss:covering_theory} and~\ref{ss:group_actions} above, on coverings of complexes of groups
and group actions on complexes of groups, and those of Section~\ref{s:unfolding} above, on unfoldings.
After proving the Density Theorem, in Section~\ref{ss:new_lattices} below we sketch how these techniques
may be used to construct uniform lattices in addition to those needed for the proof.

\subsection{Proof of the Density Theorem}

Let $X$ be a regular right-angled building of type $(W,S)$ with parameters $\{q_s\}$ (see Section~\ref{ss:rabs}).  Let
$G=\Aut(X)$ and let $\G_0 \leq G$ be the standard uniform lattice (see Section~\ref{ss:complexes_of_groups}).  Let $Y_n$
be the combinatorial ball in $X$ of radius $n \geq 0$, and let $x_0$ be the center of the chamber
$Y_0$. We first establish the following reduction:

\begin{lemma}\label{l:reduction1} To prove the Density Theorem, it suffices to show that for any $g \in \Stab_{G}(x_0)$,
and for any integer $n \geq 0$, there is a $\gamma =\gamma_n \in \CommG(\G_0)$ such that \[ g |_{Y_n} = \gamma|_{Y_n}.
\] \end{lemma}

\begin{proof}  Let $G_X(Y_0)$ be the complex of groups defined in Section~\ref{ss:rabs} above, with
fundamental group $\G_0$.  Let $G_0 = \Aut_0(X)$ be the group of type-preserving automorphisms of $X$.  Then
$G_0 \bs X$ is the chamber $Y_0$.  With the piecewise Euclidean metric on $X$ provided by
Theorem~\ref{t:metrize} above, the action of the full automorphism group $G$ on $X$ must preserve the
cardinality of types of faces in $X$.  Hence the quotient $G \bs X$ is a further quotient of $Y_0$, by the
action of the finite (possibly trivial) group of permutations  \[ H := \{ \varphi \in \Sym(S) \mid
q_{\varphi(s)} = q_s \mbox{ and } m_{\varphi(s)\varphi(t)} = m_{st} \mbox{ for all $s,t \in S$}\}.  \]

Let $Z_0 = H \bs Y_0 = G \bs X$.  By construction of $G_X(Y_0)$, the action of $H$ on $Y_0$ naturally extends
to an action by simple morphisms on the complex of groups $G_X(Y_0)$.  Let $H(Z_0)$ be the complex of groups
induced by the $H$--action on $G_X(Y_0)$.  Let $\G_0'$ be the fundamental group of $H(Z_0)$.  By
Theorem~\ref{t:group_action} above, there is an induced finite-sheeted covering of complexes of groups
$G_X(Y_0) \to H(Z_0)$.  Hence by covering theory for complexes of groups, $\G_0'$ is a finite index subgroup
of $\G_0$.  In particular, $\G_0'$ is commensurable to $\G_0$.  So $\CommG(\G_0') = \CommG(\G_0)$.

Since $G \bs X = Z_0 = \G_0' \bs X$, it follows that \[ G \bs X = \CommG(\G_0') \bs X = \CommG(\G_0) \bs X.
\]  Thus we have equality of orbits $G \cdot x_0 = \CommG(\G_0) \cdot x_0$, and so \[ G = \CommG(\G_0) \cdot
\Stab_{G}(x_0). \] Hence to show that $\CommG(\G_0)$ is dense in $G$, it is enough to show that $\CommG(\G_0)
\cap \Stab_{G}(x_0)$ is dense in $\Stab_{G}(x_0)$.  And for this, it suffices to prove the statement of this
lemma.\end{proof}

To continue with the proof of the Density Theorem, fix $g \in \Stab_{G}(x_0)$ and $n \geq 0$, and let $G_X(Y_n)$ be the canonical
complex of groups over the combinatorial ball $Y_n$, as defined in Section~\ref{ss:clumps} above.  Let $H_n$ be the finite
group obtained by restricting the action of $\Stab_G(x_0)$ on $X$ to $Y_n$, and note that $g|_{Y_n} \in H_n$.

\begin{proposition}\label{p:action_extends}   The action of $H_n$ on $Y_n$ extends to an action by simple
morphisms on the complex of groups $G_X(Y_n)$.\end{proposition}

\begin{proof}  We first show:

\begin{lemma}\label{l:sides}  For all sides $\cK$ of $Y_n$ and all $h \in H_n$, the image $h.\cK$ is also a side of
$Y_n$. That is, the action of $H_n$ takes sides to sides.  \end{lemma}

\begin{proof} Let $\cK$ be a side of $Y_n$, of type $t \in S$.  The automorphism $h$ of $Y_n$ preserves the boundary
$\partial Y_n$, so for each mirror $K_t$ contained in $\cK$, the mirror $h.K_t$ is in $\partial Y_n$ as well.  Also, $h$
preserves adjacency of mirrors in $Y_n$ (recall that two mirrors are adjacent if their intersection is of type $T$ with
$|T| = 2$).  Thus it suffices to show that if two $t$--mirrors $K_t$ and $K'_t$ of $\cK$ are adjacent, then the mirrors
$h.K_t$ and $h.K'_t$ have the same type.

Let $\phi_t$ and $\phi'_t$ be the chambers of $Y_n$ containing $K_t$ and $K'_t$, respectively.  As $K_t$ and $K'_t$ are
adjacent and of the same type, there is a unique $s \in S$, with $m_{st} = 2$, such that $\phi_t$ is $s$--adjacent to
$\phi'_t$.  Thus the images $h.\phi_t$ and $h.\phi'_t$ are $\tilde{s}$--adjacent, for some $\tilde{s} \in S$.  Hence
$(h.\phi_t,h.\phi'_t)$ is a gallery of type $\tilde{s}$ in $Y_n$.

Suppose the type of $h.K_t$ is $u$ and that of $h.K'_t$ is $u'$, with $u \neq u'$.  Since the mirrors $h.K_t$ and
$h.K'_t$ are adjacent and of distinct types, there is a chamber $\phi$ of $X$ (not necessarily in $Y_n$) which contains
both $h.K_t$ and $h.K'_t$. Thus there is a gallery $(h.\phi_t, \phi, h.\phi'_t)$ of type $(u,u')$ in $X$.  But by the
definition of the $W$--distance function on $X$, this means $\tilde{s} = uu'$, which is impossible.  Hence $u = u'$, as
required. \end{proof}

We now, for each $h \in H_n$, define a simple isomorphism of complexes of groups $\Phi^h=(\phi^h_\s):G_X(Y_n)
\to G_X(Y_n)$.  For each $s \in S$, fix a generator $g_s$ of the cyclic group $G_s = \Z/q_s\Z$.  Let $\s$ be
a vertex of $Y_n$.  By definition of the complex of groups $G_X(Y_n)$, if $\s$ is in $Y_n - \partial Y_n$
then $G_\s(Y_n)$ is the trivial group.  Now the vertex $h.\s$ is in the boundary $\partial Y_n$ if and only if
$\s \in \partial Y_n$, so for all $\s \in Y_n - \partial Y_n$ we may define the local map $\phi^h_\s:G_\s(Y_n) \to
G_{h.\s}(Y_n)$ to be the trivial isomorphism.

If $\s$ is in $\partial Y_n$ then  \[G_\s(Y_n) = G_{T_{\partial Y_n}} = \prod_{t
\in T_{\partial Y_n}} \langle g_t \rangle.\] To define the local map $\phi^h_\s$ for $\s \in \partial Y_n$,
let $T_{\partial Y_n}$ be the boundary type of $\s$ in $Y_n$ and let $U_{\partial Y_n}$ be the boundary type
of $h.\s$ in $Y_n$.  Let $t \in T_{\partial Y_n}$.  Then $\s$ is contained in a side $\cK$ of $Y_n$ of type
$t$.  By Lemma~\ref{l:sides} above, the image $h.\cK$ is a side of $Y_n$. Denote by $u_t$ the type of the
side $h.\cK$.  Since $h$ is an automorphism of $Y_n$, the map $t \mapsto u_t$ is a bijection $T_{\partial
Y_n} \to U_{\partial Y_n}$.  Since $h$ is the restriction of an automorphism of $X$ to $Y_n$,  for all $t \in
T_{\partial Y_n}$, we have $q_t = q_{u_t}$.  We may thus define the local map $\phi^h_\s: G_\s(Y_n) \to
G_{h.\s}(Y_n)$ to be the isomorphism of groups $G_{T_{\partial Y_n}} \to G_{U_{\partial Y_n}}$ induced by
$g_t \mapsto g_{u_t}$ for each $t \in T_{\partial Y_n}$.

Recall that all monomorphisms $\psi_a$ along edges in the complex of groups $G_X(Y_n)$ are the identity or natural
inclusions.  Using this, it is not hard to verify that $\Phi^h$ so defined is a simple morphism of complexes of groups.
Since $h$ is an isomorphism of $Y_n$ and each local map $\phi_\s^h$ an isomorphism of groups, it follows that $\Phi^h$
is a simple isomorphism of the complex of groups $G_X(Y_n)$.  Moreover, for all $h, h' \in H_n$, from the definition of
composition of simple morphisms (see~\cite{BH}) it is immediate that $\Phi^h \circ \Phi^{h'} = \Phi^{hh'}$.  Hence the
group $H_n$ acts on the complex of groups $G_X(Y_n)$ by simple morphisms. \end{proof}

To finish proving the Density Theorem, let $\G_n$ be the fundamental group of $G_X(Y_n)$.  By Corollary~\ref{c:Gamma_n}
above, $\G_n$ is a finite index subgroup of $\Gamma_0$.  Let $Z_n = H_n \bs Y_n$, let $H(Z_n)$ be the complex of groups
induced by the action of $H_n$ on $G_X(Y_n)$, and let $\G_n'$ be the fundamental group of $H(Z_n)$.   Since the induced
covering of complexes of groups $G_X(Y_n) \to H(Z_n)$ is finite-sheeted, $\G_n$ is a finite index subgroup of $\G_n'$.
Therefore $\G_n'$ and $\G_0$ are commensurable.  Now the group $H_n$ is, by definition, the restriction of the group
$\Stab_G(x_0)$ to the combinatorial ball $Y_n$.  Hence $H_n$ fixes the basepoint $x_0$, and so by
Theorem~\ref{t:group_action} above, $H_n$ injects into $\G_n'$.  Since $g|_{Y_n} \in H_n$, it follows that there is an element $\gamma \in \G_n'$ such that  $\gamma|_{Y_n} = g|_{Y_n}$.  But since $\G_n'$ and $\G_0$
are commensurable, $\gamma \in \CommG(\G_0)$.  By the reduction established in Lemma~\ref{l:reduction1} above, this completes the proof of the Density
Theorem.

\subsection{Constructing other lattices using unfoldings and group actions on complexes of
groups}\label{ss:new_lattices}

Let $X$ be a regular right-angled building and let $G=\Aut(X)$.  In this section we sketch how the
techniques of unfoldings and group actions on complexes of groups may be combined to construct uniform
lattices in $G$ in addition to the sequences $\G_n$ and $\G_n'$ above.

Let $Y$ be any subcomplex of $X$ obtained by unfolding the chamber $Y_0$ finitely many times.  Let $G(Y)
= G_X(Y)$ be the canonical complex of groups over $Y$ defined in Section~\ref{ss:clumps} above.  By
Proposition~\ref{p:unfolding_admissible} above, the fundamental group $\G$ of $G(Y)$ is a uniform
lattice in $G$.  The possible fundamental domains $Y$ for $\G$ include many subcomplexes which are not
combinatorial balls in $X$.

Now suppose $H$ is any (finite) group of automorphisms of the subcomplex $Y$.  As in
Proposition~\ref{p:action_extends} above, the action of $H$ on $Y$ extends to an action by simple
morphisms on the complex of groups $G(Y)$.  Let $\G'$ be the fundamental group of the induced complex of
groups over $H \bs Y$.  Then $\G'$ is also a uniform lattice in $G$.  This construction thus yields many
additional uniform lattices in $G$.

\section{Further applications of unfoldings}\label{s:applications}

In this section we give two further applications of the technique of unfoldings, which was developed in
Section~\ref{s:unfolding} above. Let $X$ be a regular right-angled building of type $(W,S)$ and
parameters $\{q_s\}$, as defined in Section~\ref{ss:rabs} above.  Let $G=\Aut(X)$ and let $G_0 =
\Aut_0(X)$ be the group of type-preserving automorphisms of $X$.  In Section~\ref{ss:discreteness} we
determine exactly when $G$ and $G_0$ are nondiscrete groups (Theorem~\ref{t:autom_gp_intro} of the
introduction).  We then
in Section~\ref{ss:strongly_transitive} prove Theorem~\ref{t:transitive} of the introduction, which
states that $G$ acts strongly transitively on $X$.

\subsection{Discreteness and nondiscreteness of $G$ and $G_0$}\label{ss:discreteness}

Let $L$ be a polyhedral complex.  Recall that $L$ is \emph{rigid} if
for any $g \in \Aut(L)$, if $g$ fixes the star in $L$ of a vertex $\s \in V(L)$, then $g=\Id_L$.  If $L$
is not rigid it is said to be \emph{flexible}.  For example, a complete graph is rigid, while a complete
bipartite graph $L=K_{q,q}$, with $q > 2$, is flexible.

The following statement is equivalent to
Theorem~\ref{t:autom_gp_intro} above.

\begin{theorem}\label{t:autom_gp} Let $X$ be a regular right-angled building of type $(W,S)$ and
parameters $\{ q_s\}$. Let $G=\Aut(X)$ and let $G_0 = \Aut_0(X)$ be the group of type-preserving
automorphisms of $X$.  Suppose $W$ is infinite and let $L$ be the nerve of $(W,S)$. \begin{enumerate}
\item\label{i:both_nondiscrete} If there are $s,t \in S$ such that $q_s > 2$ and $m_{st} = \infty$ then
$G_0$ and $G$ are both nondiscrete.  \item\label{i:mixed} If all $q_s = 2$, then $G_0$ is discrete, and
$G$ is nondiscrete if and only if $L$ is flexible. \item\label{i:last_case} If there is some $q_t > 2$,
and for all $t\in S$ with $q_{t} > 2$ we have $m_{st} = 2$ for all $s \in S - \{t\}$, then $G_0$ is
discrete, and $G$ is nondiscrete if and only if $L$ is flexible. \end{enumerate} \end{theorem}

\noindent Note that if the Coxeter group $W$ is finite then the building $X$ is finite, so both $G$ and $G_0$ are finite groups.

\begin{proof} Several results of~\cite{T1} imply that in Case~\eqref{i:both_nondiscrete}, the group $G_0$ is
nondiscrete.  For example, the set of covolumes of lattices in $G_0$ contains arbitrarily small elements.  Since a
subgroup of a discrete group is discrete, the full automorphism group $G$ is thus nondiscrete as well.

Suppose next that all $q_s = 2$.  Then $X$ is just the Davis complex $\Sigma$ for $(W,S)$. Assume $g_0 \in
G_0$ fixes a chamber $\phi$ of $X$ pointwise.  Then for each $s \in S$, since $q_s = 2$ there is a unique
chamber $\phi_s$ of $X$ such that $\phi_s$ is $s$--adjacent to $\phi$.  Since $g_0$ preserves types and fixes
$\phi$ pointwise, the element $g_0$ fixes each adjacent chamber $\phi_s$ pointwise as well.  By induction,
$g_0$ fixes the building $X$ pointwise. Hence $G_0$ is discrete.  Haglund--Paulin~\cite{HP1} and White~\cite{W} proved that the
full automorphism group $G=\Aut(X)=\Aut(\Sigma)$ is nondiscrete exactly when the nerve $L$ of $(W,S)$ is
flexible.

Suppose finally that we are in Case~\eqref{i:last_case}.  Then in particular the set $T:=\{ t \in S \mid q_t > 2 \}$ is
a nonempty spherical subset of $S$.  Let $\cC$ be the clump obtained by unfolding the
chamber $Y_0$ along all of its mirrors of types $t \in T$ (in some order).  More precisely, $\cC$ is the clump
obtained by unfolding $Y_0$ along some sequence of (possibly extended) sides of types $t \in T$, as in the proof of
Lemma~\ref{l:unfoldingYn} above.  By Proposition~\ref{p:unfolding_admissible} above, the complex of groups $G_X(\cC)$ is
admissible.  Hence $\cC$ is a strict fundamental domain for the action of a uniform lattice $\G := \pi_1(G_X(\cC))$ on
$X$, and so we may think of $X$ as tesselated by copies of $\cC$.

By Lemma~\ref{l:mirrors} above, since $T$ is a nonempty spherical subset of $S$, the union of mirrors $\cup_{t \in T}
K_t$ of $Y_0$ is contractible and thus connected.  Therefore, every mirror of $\cC$ of type $t \in T$ is in the interior
of $\cC$.  Thus every side of $\cC$ is of type $s \in S - T$.

Now suppose $g_0 \in G_0$ fixes $\cC$ pointwise.  Let $\phi$ be a chamber of $X$ which is $s$--adjacent to a chamber in
$\cC$, for some $s \in S - T$.  Then since $q_s = 2$ and $g_0$ is type-preserving, $g_0$ must fix the chamber $\phi$
pointwise.  For each $t \in T$, let $K_{\phi,t}$ be the $t$--mirror of $\phi$.  By hypothesis, $m_{st} = 2$, so the
mirror $K_{\phi,t}$ of $\phi$ is adjacent to a mirror (of type $s$) in $\partial\cC$.  Thus any chamber of $X$ which is
$t$--adjacent to $\phi$ is $s$--adjacent to a chamber in $\cC$.  Since $q_s = 2$, it follows that any chamber of $X$
which is $t$--adjacent to $\phi$, for $t \in T$, must also be fixed pointwise by the element $g_0$.  Hence $g_0$ fixes
pointwise the copy of $\cC$ in $X$ which contains the chamber $\phi$.

We have shown that for all $s \in S - T$, every copy of $\cC$ in $X$ which is $s$--adjacent to the original clump $\cC$
is also fixed pointwise by $g_0$.  By induction, $g_0 = \Id_X$.  Thus the group $G_0$ of type-preserving automorphisms
of $X$ is discrete.  The proof that $G = \Aut(X)$ is nondiscrete if and only if $L$ is flexible is by similar arguments to those of
Haglund--Paulin~\cite{HP1}. \end{proof}

\subsection{Strong transitivity}\label{ss:strongly_transitive}

We conclude by proving Theorem~\ref{t:transitive} of the introduction.  We will actually show:

\begin{theorem}\label{t:transitive_proof}  Let $X$ be a regular right-angled building of type $(W,S)$ and parameters
$\{q_s\}$, and let $G_0=\Aut_0(X)$.  Let $x_0$ be the center of the chamber $Y_0$. \begin{enumerate} \item\label{i:stab_transitive} The group $H_0 :=
\Stab_{G_0}(x_0)$ acts transitively on the set of apartments containing $Y_0$. \item The group $G_0$ acts
transitively on the set of pairs \[\{(\phi,\Sigma) \mid \Sigma \mbox{ is an apartment of $X$ containing the chamber
$\phi$}\}.\]\end{enumerate} \end{theorem}

\begin{corollary} The group $G$ acts strongly transitively on $X$. \end{corollary}

\begin{proof}[Proof of Theorem~\ref{t:transitive_proof}] Since $G_0$ acts transitively on the set of chambers
of $X$, it is enough to show~\eqref{i:stab_transitive}.   We fix an increasing sequence of subcomplexes $\cC_n$
of $X$ such that $\cC_n$ is a clump obtained by $n$ unfoldings of $\cC_0 = Y_0$ and $X = \cup_{n=0}^\infty
\cC_n$.

\begin{lemma}\label{l:transitive} Let $\Sigma$ and $\Sigma'$ be distinct apartments of $X$ which contain $Y_0$. Let $N
\geq 1$ be the smallest integer such that $\Sigma \cap \cC_{N} \neq \Sigma' \cap \cC_{N}$.  Then there is an element
$h_{N} \in H_0$ such that $h_{N}$ fixes pointwise the clump $\cC_{N-1}$, and $h_{N}(\Sigma \cap \cC_{N}) =
\Sigma' \cap \cC_{N}$.\end{lemma}

\begin{proof}  Suppose $\cC_{N}$ is obtained from $\cC_{N - 1}$ by unfolding along a side $\cK$ of type $u$.  Recall
from Lemma~\ref{l:sheets} above that $\Ch_\cK$, the set of ``new chambers" in $\cC_{N}$, consists of $q_u
- 1$ sheets.  Since $\Sigma \cap \cC_{N} \neq \Sigma' \cap \cC_{N}$, the sets of chambers $\Sigma \cap \Ch_\cK$ and
$\Sigma' \cap \Ch_\cK$ belong to different sheets in $\Ch_\cK$.

Now, for each sheet in $\Ch_\cK$, the set of chambers
in this sheet is in bijection with the set of mirrors in $\cK$.  Hence, for any two sheets in $\Ch_\cK$, there is a
type-preserving element $h'_{N} \in \Aut(\cC_{N})$ such that $h'_{N}$ fixes $\cC_{N - 1}$ pointwise, and
$h'_{N}$ exchanges these two sheets.

Since $\Sigma \cap \cC_{N - 1} = \Sigma' \cap \cC_{N - 1}$, the set of mirrors in $\cK$ contained in $\Sigma$ is
equal to the set of mirrors in $\cK$ contained in $\Sigma'$.  Thus $h'_{N}$ exchanges the sets of chambers $\Sigma
\cap \Ch_\cK$ and $\Sigma' \cap \Ch_\cK$.  So $h'_{N}$ fixes $\cC_{N - 1}$ pointwise, and $h'_{N}(\Sigma \cap
\cC_{N}) = \Sigma' \cap \cC_{N}$.

Consider the group $\langle h'_{N} \rangle $ generated by $h'_{N}$.  By similar arguments to the proof of
Proposition~\ref{p:action_extends} above, the action of $\langle h'_{N} \rangle $ on $\cC_{N}$ extends to an
action by simple morphisms on the complex of groups $G_X(\cC_{N})$.  Since the group $\langle h'_{N} \rangle$
fixes $\cC_{N-1}$ pointwise, in particular it fixes the point $x_0$.  By Theorem~\ref{t:group_action} above,
the group $\langle h'_{N} \rangle$ thus injects into the fundamental group of the induced complex of groups.
Denote by $h_{N}$ the image of $h'_{N}$ in this fundamental group.  By construction, $h_{N}$ fixes $\cC_{N -
1}$ pointwise, hence $h_{N} \in H_0$, and $h_{N}(\Sigma \cap \cC_{N}) = \Sigma' \cap \cC_{N}$ as required.
\end{proof}

Let $\Sigma$ and $\Sigma'$ be two apartments of $X$ which
contain $Y_0$.  For each $n \geq 0$ we will
construct an element $h_n \in H_0$ such that\begin{enumerate}
\item\label{i:sigmas} $h_n(\Sigma \cap \cC_n) = \Sigma' \cap \cC_n$,
and \item\label{i:clump} for all $m \geq 0$, we have
$h_{n+m}|_{\cC_{n}} = h_{n}|_{\cC_{n}}$. \end{enumerate}  Note
that, since $\cC_n \subset \cC_{n+1}$ for all $n \geq 0$, to
prove~\eqref{i:clump} it suffices to show that for all
$n \geq 0$, $h_{n+1}|_{\cC_n} = h_n|_{\cC_n}$.

To construct the sequence $\{h_n\}$, let $N  \geq 1$ be the
smallest integer such that $\Sigma \cap \cC_{N} \neq \Sigma' \cap
\cC_{N}$.  For
each $0 \leq n < N$ we define $h_n \in H_0$ to be $h_n = \Id_X$.  Let $h_{N}$
be the element of $H_0$ constructed in Lemma~\ref{l:transitive} above.
 Then for each $0 \leq n \leq N$ we have
$h_n(\Sigma \cap \cC_n) = \Sigma' \cap \cC_n$, and for all $0 \leq n <
N$ we have $h_{n+1}|_{\cC_{n}} = h_{n}|_{\cC_{n}}$.

For $n \geq N$, assume inductively that for $k \geq 0$ there are elements $h_{N}, h_{N + 1}, \ldots, h_{N +
k}$ in $H_0$ such that $h_{N + k}(\Sigma \cap \cC_{N+k}) = (\Sigma' \cap \cC_{N+k})$, and $h_{N+k}|_{\cC_{N +
k - 1}} = h_{N+k - 1}|_{\cC_{N + k - 1}}$.  To construct the next element
$h_{N + k + 1}$, note that since \[h_{N + k}(\Sigma \cap \cC_{N + k}) = \Sigma' \cap \cC_{N + k},\] the
apartments $h_{N + k}\Sigma$ and $\Sigma'$ have the same intersection with $\cC_{N + k}$.  If in addition the
apartments $h_{N + k}\Sigma$ and $\Sigma'$ have the same intersection with the next clump $\cC_{N + k + 1}$,
we put $h_{N + k + 1} = h_{N + k}$ and are done.  If not, then $N + k + 1$ is the smallest integer such that
the apartments $h_{N + k}\Sigma$ and $\Sigma'$ have distinct intersection with $\cC_{N + k + 1}$.  Hence by
Lemma~\ref{l:transitive} above, there is an element $h' \in H_0$ such that $h'$ fixes pointwise $\cC_{N +
k}$, and $h'(h_{N + k}\Sigma \cap \cC_{N + k + 1}) = \Sigma' \cap \cC_{N + k + 1}$.  We then define
$h_{N+k+1}$ to be the product $h' h_{N + k}$, and have that \[h_{N + k + 1}(\Sigma \cap \cC_{N + k + 1}) =
\Sigma' \cap \cC_{N + k + 1}.\] Since $h'$ fixes pointwise $\cC_{N + k}$, the restriction of $h_{N + k +
1}=h'h_{N + k}$ to the clump $\cC_{N + k}$ is the same as that of $h_{N + k}$.  Hence the element $h_{N + k +
1}$ has the required properties. We have thus constructed a sequence $\{h_n\}$ satisfying~\eqref{i:sigmas}
and~\eqref{i:clump} above.

By definition of the topology on $G_0$, the compact subgroup $H_0$ of
$G_0$ is complete.  The sequence $\{h_n\}$ in
$H_0$ that we have constructed is a Cauchy sequence,
by~\eqref{i:clump} above.  Hence there is an  element $h \in
H_0$ such that $h \Sigma = \Sigma'$. We conclude that $H_0$ acts
transitively on the set of apartments containing $Y_0$.
\end{proof}

\end{document}